\documentclass{amsart}
\usepackage{amsmath,amssymb}
\usepackage[all]{xy}

\newcommand{\C}{{\mathbb{C}}}

\newcommand{\N}{{\mathbb{N}}}

\newcommand{\R}{{\mathbb{R}}}

\newcommand{\Z}{{\mathbb{Z}}}

\newcommand{\Ch}{{\mathcal C}}
\newcommand{\Dh}{{\mathcal D}}

\newcommand{\Fh}{{\mathcal F}}
\newcommand{\Gh}{{\mathcal G}}

\newcommand{\Ih}{{\mathcal I}}

\newcommand{\Kh}{{\mathcal K}}

\newcommand{\Mh}{{\mathcal M}}

\newcommand{\Oh}{{\mathcal O}}

\newcommand{\Uh}{{\mathcal U}}

\newcommand{\Zh}{{\mathcal Z}}

\newcommand{\be}{\mathbf{1}}

\newcommand{\colim}{\mathrm{colim}\,}

\newcommand{\ev}{\mathrm{ev }}

\newcommand{\id}{\mathrm{id}}

\newcommand{\supp}{\mathrm{supp }\,}

\newcommand{\diag}{\mathrm{diag}}

\newcommand{\ep}{{\varepsilon}}

\newcommand{\Cs}{{$C^*$-al\-ge\-bra}}

\newcommand{\sh}{{$^*$-ho\-mo\-mor\-phism}}

\newcommand{\nprecsim}{{\operatorname{\hskip3pt \precsim\hskip-9pt |\hskip6pt}}}

\newcounter{number}[section]

\newenvironment{nummer}{\refstepcounter{number}
{\noindent\arabic{section}.\arabic{number}}}{}

\newcommand{\bn}{\begin{nummer} \rm}
\newcommand{\en}{\end{nummer}}

\newenvironment{thms}{\noindent {\sc Theorem:} \it}{}
\newenvironment{lms}{\noindent {\sc Lemma:} \it}{}
\newenvironment{props}{\noindent {\sc Proposition:} \it}{}
\newenvironment{dfs}{\noindent {\sc Definition:} \it}{}
\newenvironment{cors}{\noindent {\sc Corollary:} \it}{}

\newenvironment{rems}{\noindent {\sc Remark:}}{}

\newenvironment{exs}{\noindent {\sc Example:} }{}

\newenvironment{nproof}[1][Proof:]
{\begin{trivlist}\item[]{\sc{#1}} }
{\hbox{}\nobreak\hfill\quad\hbox{$\square$}\end{trivlist}}

\begin{document}

\title[$\Ch_{0}(X)$-algebras, stability and strongly self-absorbing $C^{*}$-algebras]
{\sc $\Ch_{0}(X)$-algebras, stability and \\
strongly self-absorbing $C^{*}$-algebras}

\author{Ilan Hirshberg}
\address{Department of Mathematics, Ben Gurion University of the
  Negev, P.O.B. 653\\
 Be'er  Sheva 84105,  Israel}

\email{ilan@math.bgu.ac.il}

\author{Mikael R{\o}rdam}
\address{Department of Mathematics and Computer Science, University of
Southern Denmark,  Campusvej 55\\
 DK-5230  Odense M, Denmark}

\email{mikael@imada.sdu.dk}

\author{Wilhelm Winter}
\address{Mathematisches Institut der Universit\"at M\"unster\\
Einsteinstr. 62\\ D-48149 M\"unster,  Germany}

\email{wwinter@math.uni-muenster.de}

\date{\today}
\thanks{{\it Research partly supported by:} Deutsche Forschungsgemeinschaft
(through the SFB 478), by the
\indent EU-Network Quantum Spaces - Noncommutative Geometry
(Contract No.\ HPRN-CT-2002-\indent00280), and by
the Center for Advanced Studies in Mathematics at Ben-Gurion
University}

\setcounter{section}{-1}

\begin{abstract}
We study permanence properties of the classes of stable and
so-called $\Dh$-stable $C^{*}$-algebras, respectively. More
precisely, we show that a $\Ch_{0}(X)$-algebra $A$ is stable if
all its fibres are, provided that the underlying compact
metrizable space $X$ has finite covering dimension or that the
Cuntz semigroup of $A$ is almost unperforated (a condition which
is automatically satisfied for $C^{*}$-algebras absorbing the
Jiang--Su algebra $\Zh$ tensorially). Furthermore, we prove that
if $\Dh$ is a $K_{1}$-injective strongly self-absorbing
$C^{*}$-algebra, then $A$ absorbs $\Dh$ tensorially if and only if
all its fibres do, again provided that $X$ is finite-dimensional.
This latter statement generalizes results of Blanchard and
Kirchberg. We also show that the condition on the dimension of $X$
cannot be dropped. Along the way, we obtain a useful
characterization of when a $C^{*}$-algebra with weakly
unperforated Cuntz semigroup is stable, which allows us to show
that stability passes to extensions  of $\Zh$-absorbing
$C^{*}$-algebras.
\end{abstract}

\maketitle

\section{Introduction}

\noindent
This paper addresses stability and $\Dh$-stability of
$\Ch_0(X)$-algebras, and related matters. Recall
(\cite{TomsWinter:ssa}) that a separable, unital, infinite
dimensional $C^*$-algebra $\Dh$ is said to be
\emph{strongly self-absorbing} if the embedding $\Dh \rightarrow
\Dh \otimes \Dh$ given by $d \mapsto d \otimes 1$ is approximately
unitarily equivalent to an isomorphism (regardless of which tensor
product is used in the definition, such a $\Dh$ is automatically nuclear,
and therefore there is no ambiguity in the
definition). The list of known examples of strongly self-absorbing
$C^*$-algebras is very short: it consists of UHF algebras of
`infinite type' (i.e., ones where all the primes which appear in
the supernatural number do so with infinite multiplicity), the
Jiang--Su algebra $\Zh$ (\cite{JiaSu:Z}), the Cuntz algebras $\Oh_2$ and
$\Oh_{\infty}$, and tensor products of $\Oh_{\infty}$ by UHF
algebras of infinite type.

We note that the algebra of compact operators on a separable
Hilbert space, $\Kh$, is not strongly self-absorbing (it is not
unital), but it does have some similar properties. In particular,
tensoring by a strongly self-absorbing $C^*$-algebra $\Dh$ can be
seen as an analogue of stabilization. We shall thus refer to a
$C^*$-algebra $A$ as being \emph{$\Dh$-stable}, or
\emph{$\Dh$-absorbing}, if $A \cong A \otimes \Dh$.
$\Dh$-stability, for the known examples of strongly self-absorbing
$C^*$-algebras, has been studied for quite some time; the concept
is of particular interest for Elliott's program to classify
nuclear $C^{*}$-algebras by their $K$-theory data (see
\cite{Ror:encyc} for an introduction). The case of UHF algebras
was studied by the second-named author in
\cite{Ror:UHF,Ror:UHFII}. Absorption of $\Oh_2$ and of
$\Oh_{\infty}$ plays an important role in the classification
results of Kirchberg and Phillips
(\cite{KirPhi:classI,Kir:fields,Phi:class}), and has been the
focus of further study (see, e.g., \cite{KirRor:pi,Kir:spi}).
Absorption of the Jiang--Su algebra seems to be particularly
important, since it is automatic for the known classifiable
algebras, whereas the known counterexamples to the Elliott
conjecture fail to absorb $\Zh$. Algebras which absorb the
Jiang--Su algebra tend to be well-behaved in many respects, and
thus, $\Zh$-stability can be seen as a regularity property, which
needs to be better understood (see \cite{GongJiangSu:Z,
Ror:Z-absorbing,TomsWinter:Zash,TomsWinter:VI,Winter:Z-class,
Winter:lfdr}). We shall comment further on this below.

Stability and $\Dh$-stability have different permanence
properties. It was shown in \cite{HjeRor:stable} that stability
for $\sigma$-unital $C^*$-algebras passes to inductive limits and
crossed products by discrete groups. Stability clearly does not in
general pass to hereditary subalgebras (look at corners of the
compact operators), and it was even shown in \cite{Ror:sns} that
there are examples of non-stable $C^*$-algebras $A$ where $M_n(A)$
is stable for some $n$. Stability passes to quotients and ideals,
however an extension of one stable $C^*$-algebra by another need
not be stable -- see \cite{Ror:extension}. We refer the reader to
\cite{Ror:stable} for a more recent survey and further details. As
for $\Dh$-stability, it was shown in \cite{TomsWinter:ssa}, under
the mild extra condition of $K_1$-injectivity (i.e., the canonical map
from $\Uh(\Dh)/\Uh_{0}(\Dh)$ to $K_{1}(\Dh)$ is injective -- a property
which holds for all known examples of  strongly self-absorbing $C^*$-algebras),
that $\Dh$-stability  passes to hereditary subalgebras, quotients,
inductive limits and extensions of separable $C^*$-algebras. Those
results have been shown earlier by Kirchberg for the cases of $\Dh
= \Oh_2$ and $\Dh = \Oh_{\infty}$, in
\cite{Kir:CentralSequences}. $\Dh$-stability does not pass in
general to crossed products; however, it does pass to crossed
products by $\Z$, $\R$ or by compact groups provided the group
action has a Rokhlin property
(\cite{HirshbergWinter:Rokhlin-ssa}).

The main questions addressed in this paper concern the behavior of
bundles of stable or $\Dh$-stable algebras.
Example~\ref{ex:nonstablebundle} (from \cite{Ror:sns}) exhibits a
non-stable $C^*$-algebra, which arises as a continuous field of
$C^*$-algebras, with fibres all isomorphic to the compacts. The
base space in that example is infinite dimensional. We shall
show, in Section~3,
that if $X$ is a locally compact Hausdorff space of finite
covering dimension, and if $A$ is a $\Ch_0(X)$-algebra such that
all its fibres are stable, then so is $A$. Thus, infinite
dimensionality is crucial in Example~\ref{ex:nonstablebundle}. We
also obtain  another regularity property of $\Zh$-stability: If
$A$ is a $\Ch_0(X)$-algebra (this time with no restriction on the
dimension of $X$), such that $A$ is $\Zh$-stable, and all the
fibres are stable, then $A$ is stable as well. To show the latter,
we use the fact (shown in \cite{Ror:Z-absorbing}) that the Cuntz
semigroup of a $\Zh$-stable $C^*$-algebra is almost unperforated,
and along the way, obtain a characterization of stability of such
$C^*$-algebras. This characterization also allows us to show that
stability is preserved under forming extensions, assuming the
algebra in the middle has almost unperforated Cuntz semigroup. We
 also use the above results to show that if $A$ is
$\Zh$-stable, and $M_n(A)$ is stable for some $n$, then so is $A$.

In Section 4,
we prove our main result that if $\Dh$ is a $K_1$-injective
strongly self-absorbing $C^*$-algebra, $A$ is a separable
$\Ch_0(X)$ algebra, $X$ is of finite covering dimension and all
the fibres of $A$ are $\Dh$-stable, then $A$ is also $\Dh$-stable.
This  generalizes earlier results of Blanchard and Kirchberg
(\cite{BlanchardKirchberg:piHausdorff}), who prove a similar
statement under the additional requirements that $\Dh =
\Oh_{\infty}$, $X = \mathrm{Prim}(A)$, and $A$ is  nuclear and
stable (with different methods; see also
\cite{BlanchardKirchberg:Glimm} for related results). The theorem
fails in the case of an infinite-dimensional base space $X$, as is
shown in Examples \ref{UHFexample} and \ref{Zexample}. In the
first example, we construct a unital $\Ch(X)$-algebra (where
$X=\prod_{\N} S^{2}$) whose fibres are isomorphic to the
UHF-algebra of type $2^\infty$ (the CAR-algebra) but which itself
does not absorb the CAR-algebra tensorially (in fact, it does not
admit a unital embedding of $M_{2}$). In the second, we give a
similar example where the fibres again are isomorphic to the
CAR-algebra, and hence are $\Zh$-stable, but where the algebra
itself is not $\Zh$-stable. On the other hand we show that if $X$
is not of finite covering dimension, but if  $A$ is \emph{locally}
$\Dh$-stable, then $A$ is also $\Dh$-stable.

In some special cases, much stronger results than ours hold: in a
very recent paper by D\u ad\u arlat, \cite{Dad:cont-fields}, it is
shown (as a corollary to a more general result) that if $A$ is a
$\Ch(X)$-algebra over a finite-dimensional compact space $X$ such
that each fibre is isomorphic to $\Dh$ -- where $\Dh$ is either
$\Oh_{2}$ or $\Oh_{\infty}$ -- then $A$ is isomorphic to $\Ch(X)
\otimes \Dh$ (the latter clearly implies that $A$ is
$\Dh$-stable). In \cite{DadWinter:asu}, D\u ad\u arlat and the third
named author derive the respective result for an arbitrary $K_{1}$-injective
strongly self-absorbing $C^{*}$-algebra.

The paper is organized as follows: in Sections 1 and 2 we recall
some facts about $\Ch_{0}(X)$-algebras and about $C^{*}$-algebras
with almost unperforated Cuntz semigroup. In Section 3 we
characterize when such algebras are stable and prove our results
about stability of $\Ch_{0}(X)$-algebras and about extensions.
Section 4 is entirely devoted to $\Dh$-stability of
$\Ch_{0}(X)$-algebras.

\section{$\Ch_{0}(X)$-Algebras}

\noindent
We recall some facts and notation about $\Ch_{0}(X)$-algebras. This is
a concept
introduced by Kasparov to generalize  continuous bundles (or fields) of
$C^{*}$-algebras over  locally compact Hausdorff spaces, cf.\
\cite{Kas:Novikov}.

\bn
\label{D-C(X)}
\begin{dfs}
Let $A$ be a $C^{*}$-algebra and $X$ a locally compact $\sigma$-compact
Hausdorff space.
$A$ is a $\Ch_{0}(X)$-algebra, if there is a  $^{*}$-homomorphism
$\mu \colon  \Ch_{0}(X) \to \Zh(\Mh(A))$ from
$\Ch_{0}(X)$ to the center of the multiplier algebra of $A$ such that,
for some (or any)
approximate unit $(h_{\nu})_{\nu \in \N}$ of $\Ch_{0}(X)$,
$\|\mu(h_{\nu})\cdot  a - a\| \to 0$ for each
$a \in A$.
\end{dfs}

The map $\mu$ is called the \emph{structure map}. We will usually
not write it explicitly. Note that if $X$ is compact, then $\mu$
has to be unital. In the compact case, we may write
``$\Ch(X)$-algebra'' instead of ``$\Ch_{0}(X)$-algebra''. \en

\bn
If $A$ is as above and  $Y \subset X$ is a closed subset, then
\[
J_{Y}:= \Ch_{0}(X \setminus Y)  \cdot A
\]
is a (closed) two-sided ideal of $A$; we denote the quotient map by
$\pi_{Y}$ and set
\[
A_{Y}:= A/J_{Y}.
\]
If $a \in A$, we sometimes write $a_{Y}$ for $\pi_{Y}(a)$. If $Y$
consists of just one point $x$, we will slightly abuse notation
and write $A_{x}$, $J_x$, $\pi_{x}$ and $a_{x}$ in place of
$A_{\{x\}}$, $J_{\{x\}}$, $\pi_{\{x\}}$ and $a_{\{x\}}$,
respectively. We say that $A_{x}$ is the \emph{fibre of $A$ at
$x$}. \en

\bn \label{upper-semicontinuity} For any $a \in A$ and for any $f
\in \Ch_0(X)$ we have $\pi_x(f \cdot a) = f(x) \pi_x(a)$ (because
$(f h_\nu - f(x)h_\nu)a$ belongs to $J_x$ for all $h_\nu$ in an
approximate unit for $\Ch_0(X)$). Moreover,
\begin{equation*}
\label{C(X)-norm}
\|a_{x}\| = \inf \{\|(1-f(x))a+f \cdot a\| \mid f \in \Ch_{0}(X) \}
\end{equation*}
for every $x \in X$. The function $x \mapsto \|a_{x}\|$ from $X$ to
$\R$, being the infimum
of a family of continuous functions from $X$ to $\R$, is upper semicontinuous. That is,
for any $x_{0} \in X$ and $\varepsilon >0$, there is a neighborhood $V \subset X$
of $x_{0}$ such that $\|a_{x}\| < \|a_{x_{0}}\| + \varepsilon$ for any $x \in V$.
Equivalently,
the set $\{ x \in X \mid \| a_{x} \| < \ep \}$ is open for every $a \in A$
and every $\ep > 0$.

\en

\bn The family $\{\pi_{x} \mid x \in X\}$ is faithful, i.e., if
$a \in A$ is such that $\|a_{x}\|= 0$ for all $x \in X$, then $a
=0$.  Indeed, if $a \in A$ is a positive element, the
$C^{*}$-algebra $B:= C^{*}(\mu(\Ch_{0}(X)), a) \subset \Mh(A)$  is
commutative. Therefore, there is a character $\chi$ on $B$ such
that $\chi(a) = \|a\|$. On the other hand, since $A$ is a
$\Ch_{0}(X)$-algebra, there is $f \in \Ch_{0}(X)$ such that $\|f
\cdot a - a \|$ is small, which shows that $\chi|_{\Ch_{0}(X)}$
is nonzero, hence a character on $\Ch_{0}(X)$. But then $\chi$
annihilates $\Ch_{0}(X \setminus \{x\})$ for some $x \in X$, so
it drops to a character $\overline{\chi}$ on $\pi_{x}(C^{*}(a))$
satisfying $\chi(a) = \overline{\chi} \circ \pi_{x}(a) $.

It follows that the map
$A \to \prod_{x \in X} A_{x}$ is injective, whence
\begin{equation*}
\label{fibre-norm}
\|a\| = \sup\{ \|a_{x}\| \mid x \in X\}
\end{equation*}
for every $a \in A$.
\en

\bn
\label{multipliers-restrictions}
If $A$ is a $\Ch_{0}(X)$-algebra and $Y \subset X$ is a closed subset, then
the quotient $A_{Y}$ may be regarded as a $\Ch_{0}(Y)$-algebra by
$g \cdot a_{Y} := (\overline{g} \cdot a)_{Y}$ for $a \in A$, $g \in \Ch_{0}(Y)$ and
some lift $\overline{g} \in \Ch_0(X)$ of $g$.

If $f \in \Ch_{0}(X)$ is a function with support in $Y$, then  the closure of
$f|_{Y} \cdot A_{Y}$ is a $\Ch_{0}(Y)$-algebra which embeds isometrically
into $A$ via
\[
f|_{Y} \cdot a_{Y} \mapsto f \cdot a
\]
for $a \in A$.
\en

\bn We are indebted to E.\ Blanchard for pointing out the following
result about tensor products of $\Ch_{0}(X)$-algebras to us (see
\cite[Proposition~3.1]{Blanchard:C(X)-remarks} -- Blanchard states
the result only for compact $X$, but the  version below follows
immediately by passing to the one-point compactification). As
usual, we denote the minimal tensor product of $A$ and $B$ by $A
\otimes B$.

\label{C(X)-products}
\begin{props}
Let $X$ and $Y$ be locally compact spaces, $A$ a $\Ch_{0}(X)$-algebra
and $B$ a $\Ch_{0}(Y)$-algebra. Suppose that either $A$ or $B$ are
exact. Then, $A \otimes B$ is
a $\Ch_{0}(X \times Y)$-algebra with fibres
\[
(A \otimes B)_{(x,y)} \cong A_{x} \otimes B_{y}
\]
for $x \in X$, $y \in Y$ and structure map $\mu$ given by
\[
\mu = \mu_{A} \otimes \mu_{B} \colon \Ch_{0}(X) \otimes \Ch_{0}(Y)
\to \Mh(A \otimes B).
\]
\end{props}

\begin{rems}
Using \cite[Proposition~IV.3.4.23]{Bla:encyc}, one can replace the assumptions
on $A$ (or $B$) by asking each fibre $A_{x}$ (or $B_{x}$) of $A$ (or of $B$) to be nuclear.
\end{rems}
\en

\bn
\label{C(X)-inductive-systems}
\begin{props}
Let $(B_{n},\varphi_{n,n+1})$ and $(\Ch(Y_{n}), \gamma_{n,n+1})$ be unital
inductive systems of $C^{*}$-algebras with limits $A$ and $\Ch(X)$, respectively.
Suppose each $B_{n}$ is a $\Ch(Y_{n})$-algebra with (unital) structure maps
$\mu_{n}\colon  \Ch(Y_{n}) \to \Zh(B_{n})$ satisfying
\begin{equation}
\label{limit-compatibility}
\varphi_{n,n+1} \circ \mu_{n} = \mu_{n+1} \circ \gamma_{n,n+1}
\end{equation}
for all $n$. Then, $A$ is a $\Ch(X)$-algebra with (unital) structure map
$\mu\colon  \Ch(X) \to \Zh(A)$ satisfying
\[
\varphi_{n,\infty} \circ \mu_{n} = \mu \circ \gamma_{n,\infty}
\]
for all $n$. If $x \in X = \lim_{\leftarrow} Y_{n}$ is a point corresponding to a
sequence $\{y_{n}\}_{n \in \N} \in \prod Y_{n}$, then
\[
A_{x} \cong \lim_{\to} (B_{n})_{y_{n}}.
\]
\end{props}

\begin{nproof}
The compatibility condition \eqref{limit-compatibility} guarantees existence of the
unital $^{*}$-homomorphism $\mu\colon  \Ch(X) \to A$. Since each $\mu_{n}$ maps to
$\Zh(B_{n})$, it is clear that $\mu$ also maps to the center of $A$, whence $A$ is a
$\Ch(X)$-algebra.

Now if $x \in X$,
then $\ev_{x} \circ \gamma_{n,\infty}$ is a character on $\Ch(Y_{n})$, hence
corresponds to a point evaluation $\ev_{y_{n}}$ for some $y_{n} \in Y_{n}$.
It is clear that
\[
\gamma_{n,n+1}(\Ch_{0}(Y_{n} \setminus \{y_{n}\})) \subset \Ch_{0}(Y_{n+1}
\setminus \{y_{n+1}\})
\]
and that
\[
\Ch_{0}(X \setminus \{x\}) = \lim_{\to} \Ch_{0}(Y_{n} \setminus \{y_{n+1}\}).
\]
>From \eqref{limit-compatibility} it follows that
\[
\varphi_{n,n+1}(\mu_{n}(\Ch_{0}(Y \setminus \{y_{n}\}) )\cdot B_{n})
\subset \mu_{n+1} ( \Ch_{0}(Y_{n+1} \setminus \{y_{n+1}\})) \cdot B_{n+1},
\]
so each $\varphi_{n,n+1}$ induces a $^{*}$-homomorphism
\[
\overline{\varphi}_{n,n+1} \colon  (B_{n})_{y_{n}} \to (B_{n+1})_{y_{n+1}}.
\]
The maps $\pi_{x} \circ \varphi_{n,\infty}\colon  B_{n} \to A_{x}$ clearly induce
a $^{*}$-homomorphism
\[
\overline{\pi}_{x} \colon   \lim_{\to} (B_{n})_{y_{n}} \to A_{x}.
\]
Conversely, the maps $\pi_{y_{n}}\colon   B_{n} \to (B_{n})_{y_{n}}$ induce a
$^{*}$-homomorphism
\[
\psi\colon   A = \lim_{\to} B_{n} \to \lim_{\to} (B_{n})_{y_{n}} .
\]
Since $\mu(\Ch_{0}(X \setminus \{x\})) \cdot A = \lim_{\to}
\mu_{n}(\Ch_{0}(Y_{n} \setminus \{y_{n}\})) \cdot B_{n}$, we have
\begin{eqnarray*}
\psi(\ker \pi_{x}) & = & \psi(\mu(\Ch_{0}(X \setminus \{x\})) \cdot A) \\
& = & \psi(\lim_{\to} \mu_{n}(\Ch_{0}(Y_{n} \setminus \{y_{n}\})) \cdot B_{n}) \\
& = & \psi(\lim_{\to} \ker \pi_{y_{n}}) \\
& = & 0.
\end{eqnarray*}
Therefore, $\psi$ induces a $^{*}$-homomorphism
\[
\psi_{x}\colon  A_{x} \to \lim_{\to} (B_{n})_{y_{n}};
\]
it is routine to check that $\overline{\pi}_{x}$ and $\psi_{x}$ are mutual inverses.
\end{nproof}
\en

\begin{rems}
With a little extra effort one can prove a nonunital version of the preceding result,
replacing unitality of the connecting maps by certain nondegeneracy conditions.
\end{rems}

\bn
In Section 4 we shall have use for the following combination of Propositions
\ref{C(X)-products} and \ref{C(X)-inductive-systems}.

\begin{lms} \label{lm:inftensor}
Let $\{A_n\}_{n \in \N}$ be a sequence of exact separable unital
\Cs s such that each $A_n$ is a $\Ch(X_n)$-algebra. Then $A =
\bigotimes_{n=1}^\infty A_n$ is a $\Ch(X)$-algebra, where $X =
\prod_{n=1}^\infty X_{n}$, with fibres
$$A_x = \bigotimes_{n=1}^\infty (A_n)_{x_n},$$
for each $x=(x_1,x_2,\dots)$ in $X$.
\end{lms}

\begin{nproof}
Set $B_{n}:= \bigotimes_{k=1}^{n} A_{k}$ and define
$\varphi_{n,n+1} \colon   B_{n} \to B_{n+1}$ by
$\varphi_{n,n+1}:= \id_{B_{n}} \otimes \be_{A_{n+1}}$.
Similarly, define $Y_{n}:=\prod_{k=1}^{n} X_{k}$ and, identifying
$\Ch(Y_{n})$ with $\bigotimes_{k=1}^{n} \Ch(X_{k})$, define
$\gamma_{n,n+1}\colon   \Ch(Y_{n}) \to \Ch(Y_{n+1})$ by
$\gamma_{n,n+1}:= \id_{\Ch(Y_{n})} \otimes \be_{X_{n+1}}$.

Repeated applications of Proposition \ref{C(X)-products} show that
each $B_{n}$ is a $\Ch(Y_{n})$-algebra with structure map
$\mu_{n}\colon  \Ch(Y_{n}) \to \Zh(B_{n})$ given by
$\bigotimes_{k=1}^{n} \nu_{k}$, where the $\nu_{k}\colon
\Ch(X_{k}) \to \Zh(A_{k})$ are the structure maps of the
$\Ch(X_{k})$-algebras $A_{k}$.

The $\mu_{k}$ clearly satisfy
\[
\varphi_{n,n+1} \circ \mu_{n} = \mu_{n+1} \circ \gamma_{n,n+1}.
\]
Therefore, Proposition \ref{C(X)-inductive-systems} applies and $A$ is a
$\Ch(X)$-algebra with structure map $\mu$ compatible with the inductive
limit structure. Here, we identify $X=\prod X_{n}$ with
$\colim Y_{n} \subset \prod Y_{n}$. Now if
$x=(x_{1},x_{2}, \ldots)$ is a point in $X$, this corresponds to the
sequence $\{(x_{1},\ldots,x_{n})\}_{n\in \N} \in \prod Y_{n}$ and by
Proposition \ref{C(X)-inductive-systems} we have
\[
A_{x} = \lim_{\to} (B_{n})_{(x_{1},\ldots,x_{n})}.
\]
Using Proposition \ref{C(X)-products} we obtain
\[
A_{x} = \lim_{\to} \bigotimes_{k=1}^{n} (A_{k})_{x_{k}}
= \bigotimes_{n=1}^{\infty} (A_{n})_{x_{n}}.
\]
\end{nproof}
\en

\bn
\label{C(X)-limits}
Let $A$ be a $\Ch_{0}(X)$-algebra. Since $X$ is locally compact and $\sigma$-compact
there is an increasing sequence $\{V_{i}\}_{i\in \N}$ of open subsets of
$X$ such that each $V_{i}$ has compact closure $K_{i}$ and $X = \bigcup_{i \in \N} V_{i}$.
Each $\Ch_{0}(V_{i}) \cdot A$ is an ideal of $A_{K_{i}}$ which, at the same time, is a
$\Ch_{0}(V_{i})$-algebra (with the obvious action of $\Ch_{0}(V_{i})$). We have
$\Ch_{0}(V_{i}) \cdot A \subset \Ch_{0}(V_{i+1}) \cdot A$ for  all $i$ and
$A = \lim_{\to} \Ch_{0}(V_{i}) \cdot A$.  This in particular shows:

\begin{props}
Let $\Ih$ be a property of (separable) $C^{*}$-algebras which passes to ideals and inductive
limits. Let $X$ be a locally compact $\sigma$-compact space and suppose  that, for any
compact subset $K$ of $X$ and any (separable) $\Ch(K)$-algebra $B$ one can show that if each
$B_{x}$, $x\in K$, satisfies $\Ih$, then so does $B$. It follows that if $A$ is a (separable)
$\Ch_{0}(X)$-algebra such that each $A_{x}$, $x \in X$, satisfies
$\Ih$, then so does $A$.
\end{props}

Property $\Ih$ in the preceding proposition could for example be
stability or $\Dh$-stability
(where $\Dh$ is a $K_{1}$-injective strongly self-absorbing $C^{*}$-algebra,
cf.\ \cite[Corollary~2.3]{Ror:stable} and \cite[Corollaries~3.1 and
3.4]{TomsWinter:ssa}).
\en

\bn \label{restriction comment} Suppose $X$ is a compact subset of
a compact space $Y$, and suppose that $A$ is a $\Ch(X)$-algebra.
The restriction map $\Ch(Y) \rightarrow \Ch(X)$, composed with the
structure map from $\Ch(X)$, gives $A$ the structure of a
$\Ch(Y)$-algebra. If $x \in X \subseteq Y$, then the fibre $A_x$
is the same regardless of whether $A$ is viewed as a
$\Ch(X)$-algebra or a $\Ch(Y)$-algebra. If $x \in Y \setminus X$,
then $A_x = 0$. \en

\section{Almost unperforated Cuntz
semigroup}

\bn
We remind the reader about the ordered Cuntz semigroup $W(A)$
associated to a \Cs{} $A$ from \cite{Cuntz:dimension} (see also
\cite{Ror:Z-absorbing}). Let $M_\infty(A)^+$ denote the (disjoint)
union $\bigcup_{n=1}^\infty M_n(A)^+$. For $a \in M_n(A)^+$ and $b \in
M_m(A)^+$ set $a \oplus b = \diag(a,b) \in M_{n+m}(A)^+$,
and write $a \precsim b$ if there is a sequence $\{x_k\}$
in $M_{m,n}(A)$ such that $x_k^*bx_k \to a$. Write $a \sim b$ if $a
\precsim b$ and $b \precsim a$. Put $W(A) = M_\infty(A)^+/\! \sim$,
and let $\langle a \rangle  \in W(A)$ be the equivalence
class containing $a$. Then $W(A)$ is an
ordered abelian semigroup when equipped with the relations:
$$\langle a \rangle + \langle b \rangle = \langle a \oplus b \rangle,
\qquad \langle a \rangle \le \langle b \rangle \iff a \precsim b,
\qquad a,b \in M_\infty(A)^+.$$

Any \sh{} $\varphi \colon  A \to B$ between \Cs s $A$ and $B$ induces a
morphism $\varphi_* \colon  W(A) \to W(B)$ by $\varphi_*(\langle a \rangle)
= \langle \varphi_n(a) \rangle$, when $a$ is a positive element in $M_n(A)$
and $\varphi_n \colon  M_n(A) \to M_n(B)$ is the natural extension of
$\varphi$. It is easy to check that $\varphi_*$ is additive and order
preserving.

If $A_0$ is a hereditary sub-\Cs{} of $A$ and if $\iota \colon  A_0 \to A$
is the inclusion mapping, then $\iota_* \colon  W(A_0) \to W(A)$ is
an order-isomorphism (from $W(A_0)$ onto $\iota_*(W(A_0))$). (This is to
say that $\iota_*$ is injective, and that for $x,y \in W(A_0)$ one has
$\iota_*(x) \le \iota_*(y)$ if \emph{and only if} $x \le y$.) We can
therefore suppress $\iota_*$ and  identify $W(A_0)$ with a
sub-semigroup of $W(A)$.
\en

\bn
Recall (from \cite{Ror:Z-absorbing}) that an ordered abelian semigroup
$(W,+,\le)$ is said to be \emph{almost unperforated} if, whenever $x,y \in
W$ and $n,m \in \N$ are such that $nx \le my$ and $n > m$, one has $x \le
y$. Equivalently, $W$ is almost unperforated if $(n+1)x \le ny$ implies $x
\le y$.

We shall write $(a-\ep)_+$ for $f_\ep(a)$, when $a$ is a positive
element in $A$ and when $f_\ep(t) = \max\{t-\ep,0\}$. If $a \precsim
b$, then for every $\ep > 0$ there exists $t \in A$ such that
$(a-\ep)_+ = t^*bt$. If $a,b$ are positive elements in $A$ with
$\|a-b\| < \ep$, then $(a-\ep)_+ \precsim b$. (See eg.\
\cite{Ror:UHFII} for this.)
\en

\bn
\label{wu-ideals-quotients}
\begin{props}
Let $0 \to I \to A \to A/I \to 0$ be a short exact sequence of separable
$C^{*}$-algebras. If $W(A)$ is almost unperforated, then so are $W(I)$ and
$W(A/I)$.
\end{props}

\begin{nproof}
Following the remarks above, we identify $W(I)$ with its image in $W(A)$.
Being almost unperforated clearly passes to sub-semigroups (with the inherited
order), so $W(I)$ is almost unperforated.

We proceed to show that $W(A/I)$ is almost unperforated. Let $x,y \in
W(A/I)$ and $n \in \N$ be given such that $(n+1)x \le ny$. We must show
that $x \le y$. Let $\pi \colon  A \to A/I$ be the quotient mapping. Since
positive elements in (matrix algebras over) $A/I$ lift to positive elements
in (matrix algebras over) $A$, upon replacing $A$ with a suitable
matrix algebra over $A$, we can assume that there are positive elements
$a,b$ in $A$ such that $x = \langle \pi(a) \rangle$ and $y = \langle
\pi(b) \rangle$. Let $\ep > 0$ be given. Let $c \otimes 1_k$ denote the
$k$-fold direct sum $c \oplus c \oplus \cdots \oplus c$. Then, as $\pi(a)
\otimes 1_{n+1} \precsim \pi(b) \otimes 1_n$, there exists an element $s$
in $M_{n,n+1}(A)$ such that
$$\pi(s)^* (\pi(b) \otimes 1_n) \pi(s) = (\pi(a)-\ep)_+ \otimes 1_{n+1} =
\pi((a-\ep)_+) \otimes 1_{n+1}.$$
Put $c = s^*(b \otimes 1_n)s - (a-\ep)_+ \otimes 1_{n+1}$. Then $c$ belongs
  to $M_{n+1}(I)$ and
\begin{eqnarray*}
(a-\ep)_+ \otimes 1_{n+1} & =  & s^*(b \otimes 1_n)s - c \; \le \:
 s^*(b \otimes 1_n)s + |c| \\ & \precsim & (b \otimes 1_n) \oplus |c|.
\end{eqnarray*}
In other words, $(n+1)\langle (a-\ep)_+ \rangle \le n \langle b \rangle +
\langle |c| \rangle \le n(\langle b \rangle + \langle |c| \rangle)$, whence
$\langle (a-\ep)_+ \rangle \le \langle b \rangle + \langle |c| \rangle$ because
$W(A)$ is assumed to be almost unperforated. Applying $\pi_*$ to the last
inequality yields $\langle (\pi(a)-\ep)_+ \rangle \le \langle \pi(b)
\rangle$. Finally, as $\ep>0$ was arbitrary we get the desired inequality: $x =
\langle \pi(a) \rangle \le \langle \pi(b) \rangle = y$.
\end{nproof}
\en

\section{Stability}

\noindent
We show here that a $\Ch_{0}(X)$-algebra $A$ is stable if and only if
all its fibres $A_x$ are stable provided that
either $X$ is finite dimensional or $A$ is
$\Zh$-stable. It is not true in general that $\Ch_{0}(X)$-algebras with
stable fibres are stable,  as shown in Example~\ref{ex:nonstablebundle}
(which essentially is taken from \cite{Ror:sns}). Enroute we
give a general characterization of when a separable \Cs{} with almost
unperforated Cuntz semigroup is stable
(Theorem~\ref{stability-characterization}). Furthermore, we show that
stability passes to extensions if the algebra in the middle has almost unperforated
Cuntz semigroup.

\bn \label{stability-char}
In \cite[Theorem 2.1]{HjeRor:stable} various (algebraic)
characterizations were given of when a separable \Cs{} $A$ is stable.
One of the equivalent conditions is as follows: for any $a \in F(A)$
there is $b \in F(A)$ such that $a \perp b$ and $a \precsim b$. Here
$F(a)$ is the
set of positive elements $a$ in $A$ for which there exists a positive
element $e$ in $A$ with $a = ae = ea$.
\en

\bn \label{lm:diverse}
\begin{lms} Let $A$ be a separable \Cs.
\begin{enumerate}
\item Let $a \in F(A)$ be given. Then there is an approximate unit
  $\{e_n\}_{n=1}^\infty$ for $A$ such that $e_1a=a$ and $e_{n+1}e_n = e_n$
  for all $n$.
\item Suppose that $A$ is stable and that $\{e_n\}_{n=1}^\infty$ is an
  approximate unit for $A$ satisfying $e_{n+1}e_n = e_n$ for all $n$. Then,
  for each $k \in \N$ there exists $\ell > k$ such that $e_1 \precsim
  e_{\ell} - e_k$.
\end{enumerate}
\end{lms}

\begin{nproof} (i). Take a positive
  contraction $f \in A$ such that $fa=a$. As $\overline{(1-f)A(1-f)}$ is a
  separable \Cs, it contains a strictly positive element $d$. Notice that
  $d \perp a$. It is easy to verify that $h := f + d$ is strictly positive
  in  $A$. Note that $\varphi(h)a = a$ whenever $\varphi$ is a continuous
  function from $\R^+$ into itself with $\varphi(0)=0$ and
  $\varphi(1)=1$. Indeed, by Weierstrass' theorem it suffices to show this
  for polynomials, and by linearity it suffices to consider the case where
  $\varphi(t) = t^k$, i.e., we must show that $(f+d)^ka=a$. The expression on the
  left-hand side expands in $2^k$ terms, the first of which is $f^ka$
  (which is $a$), and the remaining $2^k-1$ terms are of the form $wdf^\ell
  a$, where $\ell \ge 0$ and where $w$ is a word in $f$ and $d$ of length
  $k-1-\ell$. As $wdf^\ell a = wda = 0$ the claim is proved.

Take now a sequence $\{\varphi_n\}_{n=1}^\infty$ of continuous functions
$\varphi_n \colon  \R^+ \to [0,1]$ such that $\varphi_n$ is zero on
$[0,1/(n+1)]$, linear on $[1/(n+1),1/n]$, and equal to 1 on
$[1/n,\infty)$. Then the sequence $\{e_n\}$ defined by $e_n = \varphi_n(h)$
has the desired properties.

(ii). As $e_1$ belongs to $\overline{(e_2-1/2)_+A(e_2-1/2)_+}$, we have $e_1
\precsim (e_2-1/2)_+$. Assume without loss of generality that $k \ge 2$.
By stability of $A$, and because $e_k$ belongs to $F(A)$, there exists a
positive element $c \in A$ such that $c \perp e_k$ and $e_2 \precsim e_k
\precsim c$. For each
$\ell > k$ put $c_\ell := e_\ell c e_\ell = (e_\ell - e_k)c(e_\ell - e_k)
\precsim e_\ell - e_k$. Then $c_\ell \to c$. Find $\eta > 0$ such that
$(e_2-1/2)_+ \precsim (c-\eta)_+$, and find $\ell > k$ such that $\|c -
c_\ell\| < \eta$. Then
$$e_1 \; \precsim \; (e_2-1/2)_+ \; \precsim \; (c-\eta)_+ \; \precsim \;
c_\ell \; \precsim \; e_\ell - e_k.$$\
\end{nproof}
\en

\bn \label{lm:stable-bundles}
\begin{lms} Let $X$ be a compact Hausdorff space and let $A$ be
a separable $\Ch(X)$-algebra. Suppose that $A_{x}$ is stable for all $x \in
X$. Then, for each $a \in F(A)$, there is a sequence $a_1,a_2, a_3,\dots $
of positive elements in $A$ such that $a, a_1, a_2,
\dots$ are pairwise orthogonal and such that $\pi_x(a) \precsim
\pi_x(a_j)$ for all $x \in X$ and for all $j$.
\end{lms}

\begin{nproof} Let $a \in F(A)$ be given. Choose an increasing approximate unit
  $\{e_n\}$ for $A$ consisting of positive contractions for which
  $e_1a=a=ae_1$ and $e_{n+1}e_n=e_n$ for
  all $n$ (cf.\ Lemma~\ref{lm:diverse}). We show that for each $k$ there is
$\ell > k$ such that
  $\pi_x(a) \precsim \pi_x(e_\ell - e_k)$ for all $x \in X$. Let $k
  \in \N$ be fixed.

By stability of $A_x$ there is $\ell_x > k$ such that $\pi_x(e_1)
\precsim \pi_x(e_{\ell_x} -e_k)$  (cf.\ Lemma~\ref{lm:diverse}).
Take $t_x \in A$ such that
$$\|\pi_x(t_x^*(e_{\ell_x} - e_k)t_x) - \pi_x(e_1)\| < 1/2.$$
Let $W_x$ be the open neighborhood
of $x$ consisting of all points $y \in X$ for which
$$\|\pi_y(t_x^*(e_{\ell_x} - e_k)t_x) - \pi_y(e_1)\| < 1/2.$$
It follows from the relation $ae_1=a$ that $a$ belongs to
$\overline{(e_{1} -1/2)_+A(e_1 -1/2)_+}$, and
in particular that $a \precsim (e_1-1/2)_+$. Hence
$$\pi_y(a) \precsim \pi_y((e_1-1/2)_+) \precsim \pi_y(e_{\ell_x} -
e_k) \precsim \pi_y(e_{\ell} - e_k)$$
 for all $y \in W_x$ and for all $\ell \ge \ell_x$.
Refine the open cover $\{W_x\}_{x \in X}$ to a finite open
cover $\{W_{x}\}_{x \in F}$. Put $\ell = \max\{\ell_x \mid x \in
F\}$. Then $\pi_x(a) \precsim \pi_x(e_\ell - e_k)$ for all $x \in X$.

We can now take $a_j$ to be $e_{k_j} -e_{k_{j-1}+1}$, where
$2 = k_0 < k_1 < k_2 < \cdots$ is a sequence of natural
numbers chosen such that $\pi_x(a) \precsim \pi_x(e_{k_j} -e_{k_{j-1}+1})$
for all $x \in X$ and all $j$.
\end{nproof}
\en

\bn
\label{fd-stable-bundles}
\begin{props}
Let $X$ be a locally compact metrizable space of finite covering
dimension. Let $A$ be a separable $\Ch_{0}(X)$-algebra. Then, $A$
is stable if and only if $A_{x}$ is stable for all $x \in X$.
\end{props}

\begin{nproof}
Since stability passes to quotients by \cite[Corollary~2.3(ii)]{Ror:stable},
each $A_{x}$ is stable if $A$ is. We show the converse.
By \ref{C(X)-limits}, it suffices to prove the assertion for compact
$X$ of finite covering dimension.

Let $a$ be an element in $F(A)$ and let $\ep >0$ be
given. We show that there is an element $t \in A$ such that $\|a-t^*t\| \le
\ep$ and $tt^* \perp a$. This will show that $A$ is stable (by
\cite[Theorem~2.1 and Proposition~2.2]{HjeRor:stable}). Denote the
dimension of the space $X$ by $n$.

It follows from Lemma~\ref{lm:stable-bundles} that there are
positive elements $a_1,a_2, \dots, a_{n+1}$ in $A$ satisfying
$\pi_x(a) \precsim \pi_x(a_j)$ for all $x \in X$ and for all $j$,
and such that $a,a_1,a_2, \dots, a_{n+1}$ are pairwise orthogonal.

For each $x \in X$ there are elements $s_{j,x}$
in $A$ such that $\|\pi_x(s_{j,x}^*a_js_{j,x} -a)\| < \ep$ for
$j=1,2, \dots, n+1$. Put $t_{j,x} :=
a_j^{1/2}s_{j,x}$. Then
$$\|\pi_x(t_{j,x}^*t_{j,x} -a)\| < \ep, \qquad t_{j,x}t_{j,x}^* \in
\overline{a_jAa_j}.$$
For each $x \in X$, let $U_x$ be the open neighborhood of $x$
consisting of all $y \in X$ for which $\|\pi_y(t_{j,x}^*t_{j,x} -a)\|<
\ep$ for all $j=1,2, \dots, n+1$.

Because $X$ has dimension $n$, the open cover $\{U_x\}_{x \in X}$ of $X$ has an
open subcover $\{V_{j,\alpha}\}$,  with $j=1,2, \dots, n+1$ and with
$\alpha$ in some
(finite) index set $I_j$, such that the sets $\{V_{j, \alpha}\}_{\alpha \in
I_j}$ are pairwise disjoint for each fixed $j$. Relabeling the elements
$t_{j,x}$ we get elements $t_{j,i,\alpha}$ in $A$ such that
$$\|\pi_y(t_{j,i,\alpha}^*t_{j,i,\alpha} -a)\| < \ep, \qquad
t_{j,i,\alpha}t_{j,i,\alpha}^* \in \overline{a_jAa_j},$$
for all $i,j=1,2, \dots, n+1$, for all $\alpha \in I_i$, and for all $y \in
V_{i,\alpha}$. Let $\{\varphi_{i,\alpha}\}$ be a partition of the unit
subordinate to the cover $\{V_{i,\alpha}\}$, i.e., $\sum_{i,\alpha}
\varphi_{i,\alpha} = 1$ and $\supp(\varphi_{i,\alpha}) \subseteq
V_{i,\alpha}$. Put
$$t := \sum_{j=1}^{n+1} \sum_{\alpha \in I_j} \varphi_{j,\alpha}^{1/2} \,
t_{j,j,\alpha}.$$
As $at_{j,j,\alpha} = 0 = t_{i,i,\beta}^*a$ for all $i,j$ and
$\alpha, \beta$, we see that $tt^* \perp a$. Next, using that
$t_{j,j,\alpha}^*t_{i,i,\beta} = 0$ if $i \ne j$ and that
$\varphi_{j,\alpha}^{1/2}\varphi_{j,\beta}^{1/2} = 0$ if $\alpha \ne
\beta$, we get
\begin{eqnarray*}
t^*t &=& \sum_{j,\alpha} \sum_{i,\beta}
\varphi_{j,\alpha}^{1/2}\varphi_{i,\beta}^{1/2} \,
t_{j,j,\alpha}^*t_{i,i,\beta} \\
&=& \sum_{j,\alpha} \varphi_{j,\alpha} \,  t_{j,j,\alpha}^*t_{j,j,\alpha}.
\end{eqnarray*}
If $\varphi_{j,\alpha}(x) \ne 0$, then $x \in V_{j,\alpha}$ and
$\|\pi_x(t_{j,j,\alpha}^*t_{j,j,\alpha} - a)\| < \ep$. Hence
$$\|\pi_x\big(\varphi_{j,\alpha} (  t_{j,j,\alpha}^*t_{j,j,\alpha}
-a)\big) \| = \varphi_{j,\alpha}(x) \|\pi_x( t_{j,j,\alpha}^*t_{j,j,\alpha}
-a)\| \le \ep \varphi_{j,\alpha}(x),$$
for all $x \in X$. We conclude that $\|\pi_x(t^*t-a)\| \le \ep$ for
all $x \in X$, so that
$$\|t^*t-a\| = \sup_{x \in X} \| \pi_x(t^*t-a)\| \le \ep$$
as desired.
\end{nproof}
\en

\bn \label{corner-of-matrix-algebra-1} It was shown in
\cite{Ror:sns} that there is a non-stable separable $C^*$-algebra
$A$ such that $M_n(A)$ is stable for some $n$ (see Example
\ref{ex:nonstablebundle} below). We shall now use Proposition
\ref{fd-stable-bundles} to show that this phenomenon cannot occur
if $A$ is $\Zh$-stable.

\begin{cors}
Suppose $A$ is a separable, $\Zh$-stable $C^*$-algebra. If
$M_n(A)$ is stable for some $n$, then so is $A$.
\end{cors}

\begin{proof}
By \cite[Proposition 2.1]{Ror:sns}, if $M_n(A)$ is stable, then
so is $M_{n+1}(A)$. Denote $$\Zh_{n,n+1} = \{f \in \Ch([0,1],M_n
\otimes M_{n+1}) \mid f(0) \in M_n \otimes 1 \, , \,  f(1) \in 1
\otimes M_{n+1}\}.$$ We know that there is a unital embedding
$\iota$ of $\Zh_{n,n+1}$ into $\Zh$ (see \cite{JiaSu:Z}). Consider
the inductive system
$$\xymatrix{
 A \otimes \Zh_{n,n+1}
\ar[r]^{id\otimes \iota} & A \otimes \Zh \ar[r]^<<<<<<{x \mapsto x
\otimes 1} & A \otimes \Zh \otimes \Zh_{n,n+1} \ar[r]^{id\otimes
\iota} & A \otimes \Zh \otimes \Zh \ar[r]^<<<<<<{x \mapsto x
\otimes 1} & \dots}
$$
The inductive limit is the same as that of
$$
\xymatrix{A \ar[r]^<<<<<<{x \mapsto x \otimes 1} & A \otimes \Zh
\ar[r]^<<<<<<{x \mapsto x \otimes 1} & A \otimes \Zh \otimes \Zh
\ar[r]^<<<<<<{x \mapsto x \otimes 1} & \dots}
$$
which is isomorphic to $A$. By skipping all the even places in
the first diagram, and noting that $A \cong A \otimes \Zh$, we
see that $A$ can be written as an inductive limit of algebras of
the form $A \otimes \Zh_{n,n+1}$. We view $A \otimes \Zh_{n,n+1}$
as a $\Ch([0,1])$-algebra, by taking the embedding of $\Ch([0,1])$
into $\Zh_{n,n+1}$ as scalar functions. The fibres are $A \otimes
M_n$ at 0, $A \otimes M_{n+1}$ at 1, and $A \otimes M_n \otimes
M_{n+1}$ elsewhere. In particular, the fibres are all stable. By
Proposition \ref{fd-stable-bundles}, $A \otimes \Zh_{n,n+1}$ is
stable. Since stability passes to inductive limits, it follows
that $A$ is stable, as required.
\end{proof}

In Corollary \ref{corner-of-matrix-algebra-2} we shall give an
alternative proof of Corollary \ref{corner-of-matrix-algebra-1},
using Theorem \ref{stability-characterization} rather than
Proposition \ref{fd-stable-bundles}.
 \en

\bn \label{stability-characterization} We give below a
characterization of stability for separable \Cs s whose Cuntz
semigroup of positive elements is almost unperforated. This result is
very similar to \cite[Proposition~5.1]{HjeRor:stable}. We refer to
\cite{BlaHan:quasitrace} for the definition and properties of quasi-traces and
2-quasi-traces.

\begin{thms} Let $A$ be a separable $C^{*}$-algebra for which $W(A)$ is
  almost unperforated. Then $A$ is stable if and only if $A$ has no
  bounded non-zero $2$-quasi-trace and $A$ has no non-zero unital
  quotient.
\end{thms}

\begin{nproof} The ``only if'' part holds without the assumption that
  $W(A)$ is almost unperforated: If $A$ is stable, then so is any
  non-zero quotient of $A$, and no (non-zero) unital \Cs{} is
  stable. Let now $\tau$ be a non-zero 2-quasi-trace on $A$. Then
  $\tau(a) > 0$ for some $a \in F(A)$. By stability there is a sequence
  $\{a_n\}_{n=1}^\infty$ of pairwise orthogonal elements in $A$ each
  equivalent to $a$, and hence with $\tau(a_n) = \tau(a)$ for all
  $n$. (Two positive elements $a$ and $b$ in a \Cs{} are \emph{equivalent} if
    $a = x^*x$ and $b = xx^*$ for some $x$ in the \Cs.)
    But then $\tau$ cannot be bounded.

To prove the ``if'' part, we show that for each $a \in F(A)$ there is
$b \in A^+$
  such that $a \perp b$ and $a \precsim b$ (cf.\ \cite[Theorem~2.1 and
  Proposition~2.2]{HjeRor:stable}). Accordingly, take $a$ in
  $F(A)$, and let $B$ be the hereditary sub-\Cs{} of $A$ consisting of
  all $x \in A$ for which $xx^* \perp a$ and $x^*x \perp a$. We find a
  positive element $b \in B$ such that $a \precsim b$.

By a \emph{state on $W(A)$ normalized at $x \in W(A)$} we mean an
additive order preserving map from $W(A)$ into $\R^+ \cup \{\infty\}$
that maps $x$ to 1. The set of all such states is denoted by $S(W,x)$.
By
  \cite[Proposition~3.2]{Ror:Z-absorbing} (essentially an argument from
  \cite{GooHan:extending})
and the assumption that $W(A)$ is almost
  unperforated, we can conclude that $a \precsim b$ (or, equivalently,
  $\langle a \rangle \le \langle b \rangle$) if $\langle a \rangle \le N
  \langle b \rangle$ for some $N$ and if $f(\langle a \rangle) < f(\langle
  b \rangle)$ for all
 $f \in S(W(A),\langle b \rangle)$.

Since $a \in F(A)$ we can find $e,e' \in F(A)$ such that $ae=a$ and
$e e'=e$.

We remark that $B$ is full in $A$. Suppose, to reach a
contradiction, that the closed two-sided ideal $I$ generated by $B$ is
proper. Then $e + I$ is a unit for $A/I$ contrary to our
assumptions. (Indeed, for all $x \in A$ we have $ex + I = x +I = xe +
I$. To see the former identity, put $y = ex-x$ and note that $yy^*
\perp a$, so $yy^*$ belongs to $B$, whence $y$ belongs to $I$.)

We next show that $B$ contains an element $b_0$ such that $\langle e
\rangle \le N \langle b_0 \rangle$ for some $N$. The set $F(A)$ is
contained in the Pedersen ideal of $A$, so $e$ belongs to the
algebraic ideal generated by $B$. Hence there exist $b_1, \dots, b_N$
in $B^+$ and $x_1, \dots, x_N$ in $A$ such that $e \le \sum_{j=1}^N
x_j^*b_jx_j$. Put $b_0 = \sum_{j=1}^N b_j$. Then $\langle b_j \rangle \le
\langle b_0 \rangle$ for all $j$, so $\langle e \rangle \le
\sum_{j=1}^N \langle b_j \rangle \le N \langle b_0 \rangle$.

We now show that
\begin{equation} \label{eq:1}
\sup \{ f(\langle b \rangle) \mid b \in B^+ \} = \infty
\end{equation}
for all $f \in S(W(A), \langle e \rangle)$. Suppose, to reach a
contradiction, that $f_0 \in S(W(A), \langle e \rangle)$ is such that
the supremum in \eqref{eq:1} is finite, say equal to $C$. Then $f_0(\langle x
\rangle) \le C+1$ for all $x \in A^+$. Indeed,
\begin{eqnarray*}
x & = & x^{1/2}ex^{1/2} +  x^{1/2}(1-e)x^{1/2} \: \precsim \;
x^{1/2}ex^{1/2} \oplus x^{1/2}(1-e)x^{1/2} \\ & \sim &
e^{1/2}xe^{1/2} \oplus (1-e)^{1/2}x(1-e)^{1/2} \; \precsim \;  e \oplus
(1-e)x(1-e),
\end{eqnarray*}
cf.\ \cite[Lemma~2.8]{KirRor:pi},
so $f_0(\langle x \rangle) \le f_0(\langle e \rangle) + C = 1+C$, because
$(1-e)x(1-e)$ belongs to $B$.

Let $\overline{f}_0$ be the lower semicontinuous dimension function
arising from $f_0$, i.e.,
$$\overline{f}_0(\langle c \rangle) = \lim_{\ep \to 0^+} f_0(\langle
(c-\ep)_+ \rangle), \qquad c \in M_\infty(A).$$
Then there is a 2-quasi-trace $\tau$ on $A$ such that
$$\overline{f}_0(\langle c \rangle) = \lim_{n \to \infty}
\tau(c^{1/n})$$
for $c \in M_\infty(A)$. To reach a contradiction with our assumptions
we show that $\tau$ is non-zero and bounded. To see the former,
note that $1 = f_0(\langle e \rangle) \le \overline{f}_0(\langle e'
\rangle) \le f_0(\langle e' \rangle) \le C+1 < \infty$. The latter follows
from the formula
$$\tau(c) = \int_0^{\|c\|} \overline{f}_0(\langle (c-t)_+ \rangle) \,
dt, \qquad c \in A^+,$$
which shows that $\|\tau\| \le C+1$. This completes the proof of
\eqref{eq:1}.

The set $S(W(A), \langle e \rangle)$ is compact when equipped with the
topology of point-wise convergence (the weak-$^*$-topology). We can
therefore find $b'_1, \dots, b'_n$ in $B^+$ such that for each $f \in
S(W(A), \langle e \rangle)$ there is at least one $j$ for which
$f(\langle b'_j \rangle) > 2$. Put $b = b_0 + b'_1 + \cdots +
b'_n$. Then $\langle a \rangle \le \langle e \rangle
\le N \langle b_0 \rangle \le N \langle
b \rangle$; and $\langle b'_j \rangle \le \langle b \rangle$ for all $j$, so
$f(\langle b \rangle) \ge 2$ for all $f \in S(W(A), \langle e
\rangle)$.

To complete the proof we must show that $f(\langle a \rangle) <
f(\langle b \rangle)$ for all $f \in S(W(A), \langle b
\rangle)$. Take such a state $f$, and note that $f(\langle e \rangle)
\le N f(\langle b \rangle) = N < \infty$. If $f(\langle e \rangle) =
0$, then $f(\langle a \rangle) = 0 < 1 = f(\langle b \rangle)$, because
$a \precsim e$. And if  $f(\langle e \rangle) = \alpha > 0$, then
$\alpha^{-1}f$ belongs to  $S(W(A), \langle e \rangle)$, in which case
we have
$$\alpha^{-1} f(\langle a \rangle) \le \alpha^{-1} f(\langle e
\rangle) = 1 < 2 < \alpha^{-1} f(\langle b \rangle),$$
as desired.
\end{nproof}
\en

\bn \label{qt-quotient} The proposition below is contained in Haagerup's
manuscript \cite{Haa:quasi}. We restate it here and give a short proof
based on a result from \cite{BlaHan:quasitrace}.

\begin{props} Any 2-quasi-trace $\tau$ defined on a \Cs{} $A$, which
  vanishes on a
  closed two-sided ideal $I$ in $A$, factors through the quotient $A/I$,
  i.e., there is a 2-quasi-trace $\overline{\tau}$ on $A/I$ such that $\tau =
  \overline{\tau} \circ \pi$, where $\pi \colon  A \to A/I$ is the quotient mapping.
\end{props}

\begin{nproof} We wish to define $\overline{\tau}$ by $\overline{\tau}(a+I) = \tau(a)$
  for $a \in A$; and we must check that this is well-defined, i.e., we must
  show that $\tau(a+x) = \tau(a)$ for all $a \in A$ and for all $x \in
  I$. Let $\{e_\alpha\}$ be an increasing approximate unit for $I$
  consisting of positive contractions, and such that $\{e_\alpha\}$ is
  asymptotically central for $A$. As $\tau$ is continuous (cf.\
  \cite[Corollary~II.2.5]{BlaHan:quasitrace}) it suffices to show that
  $\tau((1-e_\alpha)a(1-e_\alpha)) \to \tau(a)$ for all $a \in A$ (because
  $\|(1-e_\alpha)(a+x)(1-e_\alpha)- (1-e_\alpha)a(1-e_\alpha)\| \to
  0$). As quasi-traces by definition are self-adjoint, it suffices to show
  this for self-adjoint elements $a \in A$. Fixing such an element $a$, put
  $b_\alpha = (1-e_\alpha)a(1-e_\alpha)$ and $x_\alpha = a -
  b_\alpha$. Then $x_\alpha$ belongs to $I$, whence
\begin{eqnarray*}
|\tau(a) - \tau(b_\alpha)| & = & |\tau(a) - \tau(a - x_\alpha)| \\
 & \le & |\tau(a) - \tau(a-x_\alpha) - \tau(x_\alpha)| + |\tau(x_\alpha)| \\
&=&  |\tau(a) - \tau(a-x_\alpha) - \tau(x_\alpha)| \: \to \:  0,
\end{eqnarray*}
by \cite[Corollary~II.2.6]{BlaHan:quasitrace} because $\|ax_\alpha
-x_\alpha a\| \to 0$.
\end{nproof}
\en

\bn
The proposition below is due to Uffe Haagerup (private
communication).

\label{bundle-qt}
\begin{props}
Let $X$ be a compact Hausdorff space (not necessarily of finite
dimension) and let $A$ be a separable $\Ch(X)$-algebra. If $A$ admits
a bounded non-zero $2$-quasi-trace, then so does $A_{x}$ for some $x$
in $X$.
\end{props}

\begin{nproof}
Let $QT(A)$ be the compact convex set of all $2$-quasi-traces
on $A$ of norm 1, and suppose that $QT(A)$ is non-empty. Let $\tau$ be
an extreme point in $QT(A)$. We show that there is a functional
$\rho$ on $\Ch(X)$ such that $\tau(fa) = \rho(f)\tau(a)$ for all $f
\in \Ch(X)^+$ and for all $a \in A$.

To this end let $\Mh$ denote the set of positive contractions in
$\Ch(X)$, and fix an element $f$ in $\Mh$. Put $\tau_1(a)
= \tau(fa)$ and $\tau_2(a) = \tau((1-f)a)$. As $fa$ and $(1-f)a$
commute for all (self-adjoint) $a$ we see
that $\tau = \tau_1 + \tau_2$ (this identity is first verified on
self-adjoint elements, and then extended to all elements in $A$ using that
$\tau$ is self-adjoint).
The norm of a 2-quasi-trace $\sigma$ on $A$ is equal to
$\sup_n \sigma(e_n)$, where $\{e_n\}$ is any increasing approximate
unit for $A$ consisting of positive contractions. In particular,
$\|\tau\| = \|\tau_1\| + \|\tau_2\|$. By extremality of $\tau$ we
conclude that $\tau_1$ is proportional to $\tau$, and so there exists
a constant $\rho(f)$ such that $\tau_1(a) =
\rho(f)\tau(a)$ for all $a \in A$.
This shows that $\tau(fa) = \rho(f)\tau(a)$ for all $f \in
\Mh$ and for all $a \in A$.

Let $a \in A^+$ with $\tau(a) >0$ be given. Put $\rho_a(f) =
\tau(fa)/\tau(a)$. Then $\rho_a \colon  \Ch(X) \to \C$ is linear and
$\rho_a(f) = \rho(f)$ for all $f \in \Mh$. As a linear functional on
$\Ch(X)$ is determined by its values on $\Mh$ it follows that $\rho_a$
is independent of the choice of $a \in A^+$, and that $\rho$ extends
uniquely to a linear functional on $\Ch(X)$, again denoted by
$\rho$, such that $\tau(fa) = \rho(f)\tau(a)$ for all $f \in \Ch(X)$ and all
$a \in A^+$ with $\tau(a) > 0$. If $a \in A^+$ and $\tau(a) = 0$,
then $\tau(fa) = 0 = \rho(f)\tau(a)$ for all $f \in \Ch(X)$. In conclusion
we have $\tau(fa) = \rho(f)\tau(a)$ whenever $f \in \Ch(X)$ and $a \in
A^+$, or (since $\tau$ is self-adjoint and additive on commuting
self-adjoint elements) for all $f \in \Ch(X)^+$ and for all $a \in A$.

Fix a positive element $a$ with $\tau(a) >0$, and let $f,g \in
\Ch(X)$. Then
$$\rho(fg)\tau(a) = \tau(fga) = \rho(f)\tau(ga) =
\rho(f)\rho(g)\tau(a).$$
This shows that $\rho$ is multiplicative.
It follows that there exists $x \in X$ such that $\rho(f) = f(x)$ for all
$f \in \Ch(X)$.

We can now conclude that $\tau$ vanishes on $I_{x} = \Ch_0(X \setminus
\{x\}) \cdot A$. It follows from
Proposition~\ref{qt-quotient}  that
$\tau$ drops to the quotient $A_x$, i.e., that there is a
2-quasi-trace $\overline{\tau}$ on $A_x$ such that $\tau = \overline{\tau} \circ
\pi_x$. In particular, $\overline{\tau}$ is a bounded non-zero
2-quasi-trace on $A_x$.
\end{nproof}
\en

\bn
\label{stable-bundles}
\begin{props}
Let $X$ be a locally compact Hausdorff space (not necessarily of finite
dimension) and let $A$ be a separable $\Ch_{0}(X)$-algebra. Suppose that
$W(A)$ is almost unperforated (for example, $A$ could be $\Zh$-stable). Then,
$A$ is stable if and only if $A_x$ is stable for each $x \in X$.
\end{props}

\begin{nproof}
As in the proof of Proposition \ref{fd-stable-bundles}, it is clear that
stability passes from $A$ to each fibre $A_{x}$. Again by \ref{C(X)-limits}
(and using Proposition \ref{wu-ideals-quotients}),
it suffices to prove the converse for compact $X$.

By Theorem~\ref{stability-characterization} it suffices to show that
$A$ has no non-zero bounded 2-quasi-trace and that $A$ has no non-zero
unital quotient. The former follows from Proposition~\ref{bundle-qt} since no
fibre $A_x$ admits a non-zero bounded 2-quasi-trace (again by
Theorem~\ref{stability-characterization}).

Suppose that $J$ is a
closed two-sided ideal in $A$ such that $A/J$ is unital. Then
$A_x/\pi_x(J)$ is unital and at the same time stable (being a quotient of $A/J$,
cf.\ \cite[Corollary~2.3(ii)]{Ror:stable}),
and is therefore necessarily zero (by
Theorem~\ref{stability-characterization}). This proves that $\pi_x(J) = A_x$ for all
$x$, whence $J = A$. To see the latter, note that if $J \ne A$, then
there is a positive contraction $a \in
  A$ such that $\|a + J\| =1$. Let $\{e_n\}$ be an approximate unit
  for $J$. Then $\|(1-e_n)a(1-e_n)\| = 1$ for all $n$. So there is
  $x_n \in X$ with $\|\pi_{x_n}((1-e_n)a(1-e_n))\| \ge 1/2$. Let $x_0$
be an accumulation point for $\{x_n\}$. Then
$\|\pi_{x_0}((1-e_n)a(1-e_n))\| \ge 1/2$ for all $n$. But then
$\pi_{x_0}(a) \notin \pi_{x_0}(J)$.
\end{nproof}
\en

\bn
\label{stable-extensions}
\begin{cors}
Let $0 \to J \to A \to B \to 0$ be a short exact sequence of separable
$C^{*}$-algebras.
\begin{enumerate}
\item If $J$ and $B$ are both stable and  $W(A)$ is almost
  unperforated (for example, $A$ could be $\Zh$-stable), then $A$ is
  stable.
\item If both $J$ and $B$ are stable and $\Zh$-stable, then $A$ is
  also stable and $\Zh$-stable.
\end{enumerate}
\end{cors}

\begin{nproof}
(i). In view of Theorem \ref{stability-characterization} we have to
show that $A$ has no bounded non-zero $2$-quasi-trace and no non-zero
unital quotient. Let $\pi \colon  A \to B$ denote the quotient mapping.

Suppose, to reach a contradiction, that $I$ is a proper closed
two-sided ideal in $A$ and that $A/I$ is
unital. Then we have the following commutative diagram with exact
rows:
$$\xymatrix{0 \ar[r] & J \ar[r] \ar[d] & A \ar[r]^-\pi \ar[d] & B
  \ar[r] \ar[d] & 0 \\  0 \ar[r] & J/(J \cap I) \ar[r]  & A/I \ar[r] & B/\pi(I)
  \ar[r] & 0.}$$
We have either $\pi(I) = B$ or $B/\pi(I)$ is unital (and
non-zero). The latter is impossible because $B$ is stable. Hence
$\pi(I) = B$, in which case $J/(J \cap I)$ is isomorphic to $A/I$. The
latter is unital and non-zero; hence $J/(J \cap I)$ is unital and
non-zero, contradicting that $J$ is stable.

Suppose next that $\tau$ is a bounded $2$-quasi-trace on $A$. The
restriction of $\tau$ to $J$ is also a bounded $2$-quasitrace, and is hence
zero by Theorem~\ref{stability-characterization}. By
Proposition~\ref{qt-quotient},  $\tau$ drops
to a $2$-quasi-trace $\overline{\tau}$ on $B$. Again by
Theorem~\ref{stability-characterization} we conclude that $\overline{\tau}
= 0$. This entails that $\tau = \overline{\tau} \circ \pi = 0$.

As $A$ has no (non-trivial) unital quotient and no non-zero
bounded 2-quasi-trace, Theorem~\ref{stability-characterization}
yields that $A$ is stable.

(ii). By \cite[Theorem~4.3]{TomsWinter:ssa}, $\Zh$-stability passes to extensions,
whence $A$ is $\Zh$-stable. Now, $W(A)$ is almost unperforated by
\cite[Theorem~4.5]{Ror:Z-absorbing}, and (i) applies.
\end{nproof}

It should be noted that not all extensions of stable (separable) \Cs s
are stable, cf.\ \cite{Ror:extension}.
\en

\bn
\begin{exs} \label{ex:nonstablebundle} We mention here an example from
\cite[Section~4]{Ror:sns} of a non-stable $\Ch(X)$-algebra $B$ whose
fibres $B_x$ are isomorphic to $\Kh$ for all $x \in X$. This example shows
that Propositions~\ref{fd-stable-bundles} and \ref{stable-bundles}
cannot be improved by removing the condition that $X$ is of finite
dimension in the former or that the algebra is $\Zh$-stable in the latter.

In the example $X$ is an infinite cartesian product of Moore
spaces $Y_n$, where $n$ can be taken to be any integer $\ge 2$.
The \Cs{} $B$ is the hereditary sub-\Cs{} of $\Ch(X,\Kh)$,
$$B = \overline{\bigcup_{n=1}^\infty (p_1 \oplus p_2 \oplus \cdots \oplus
  p_n)\Ch(X,\Kh)(p_1 \oplus p_2 \oplus \cdots \oplus p_n)},$$
where $\{p_j\}$ is a certain sequence of 1-dimensional projections in
$\Ch(X,\Kh)$. The claim to fame of $B$ in the context of
\cite{Ror:sns} is that $M_k(B)$  is non-stable for $1 \le k < n$, but
$M_n(B)$ is stable.

Any hereditary sub-\Cs{} $B$ of $\Ch(X,\Kh)$ is a
$\Ch(X)$-algebra. The fibre map $\pi_x \colon  B \to B_x$ coincides
with the restriction of the evaluation mapping at $x$ to $B$. Hence
$B_x = \pi_x(B)
\subseteq \Kh$. In the case at hand, $B_x \cong \Kh$ for all $x \in
X$ (because $B_x$ is an infinite dimensional hereditary sub-\Cs{} of
$\Kh$). Hence all fibres of $B$ are stable, but $B$ itself is not stable.
\end{exs}
\en

\bn \label{corner-of-matrix-algebra-2}
\begin{cors}
Suppose $A$ is a separable $C^*$-algebra such $W(A)$ is almost
  unperforated (for example, $A$ could be $\Zh$-stable). If $M_n(A)$ is
stable for some $n$, then so is $A$.
\end{cors}

\begin{nproof}
By \cite[Proposition II.4.1]{BlaHan:quasitrace}, a bounded nonzero
2-quasi-trace on $A$ extends to one on $M_n(A)$. If $A$ had a
non-zero unital quotient, then clearly so would $M_n(A)$. Thus,
since $M_n(A)$ is stable, it follows that $A$ has no unital
quotients and does not admit a bounded, non-zero 2-quasi-trace. By
Theorem \ref{stability-characterization}, $A$ is indeed stable.
\end{nproof}
 \en

\section{$\Dh$-stability}

\noindent In this section we show that, for a $K_{1}$-injective
strongly self-absorbing $C^{*}$-algebra $\Dh$, a
$\Ch_{0}(X)$-algebra $A$ is $\Dh$-stable if and only if all its
fibres  are, provided that $X$ is finite-dimensional (Theorem
\ref{D-stable-bundles}). We provide examples (\ref{UHFexample} and
\ref{Zexample}) showing that the above statement can fail with $X$
infinite dimensional, for $\Dh$ a UHF algebra or the Jiang--Su
algebra. However, we show in Proposition \ref{locally D-stable}
that if $A$ is `locally' $\Dh$-stable then $A$ must be
$\Dh$-stable (even when $X$ is infinite dimensional).

\bn \label{D-stability-characterization} To each \Cs{} $A$ one
associates the \Cs{} $\prod_\N A$ of all bounded sequences in $A$,
the \Cs{} $\bigoplus_\N A$ of all sequences in $A$ that converge
to zero, and the quotient $A_{\infty} = \prod_{\N}
A/\bigoplus_{\N} A$. We view $A$ as embedded in $A_{\infty}$ as
the (equivalence classes of) constant sequences.

\begin{props}
Let $A$ and $\Dh$ be separable $C^{*}$-algebras, such that $\Dh$
is $K_{1}$-injective and strongly self-absorbing. Then, the
following are equivalent:
\begin{enumerate}
\item[a)] $A$ is $\Dh$-stable.
\item[b)] Given $\eta>0$ and finite subsets $\Fh \subset A$ and
$\Gh \subset \Dh$, there is a c.p.c.\ map $\psi\colon \Dh \to A$ such that
\begin{enumerate}
\item[(i)] $\|b \psi(\be_{\Dh}) - b\| < \eta$
\item[(ii)] $\|b \psi(d) - \psi(d) b\| < \eta $
\item[(iii)] $\| b(\psi(dd') - \psi(d) \psi(d'))\|< \eta$
\end{enumerate}
for all $b \in \Fh, \, d,d' \in \Gh$.
\item[c)] Given $\eta>0$ and a finite subset $\Fh \subset A$,
there are a $^{*}$-homomorphism $\kappa\colon A \to A$ and a unital
$^{*}$-homomorphism $\mu\colon \Dh \to \Mh(A)$ (the multiplier algebra of $A$)
such that
\[
 [\kappa(A),\mu(\Dh)]=0  \mbox{ and } \|\kappa(b)-b\| < \eta
\]
for all  $b \in \Fh$.
\item[d)] There exists a c.p.c.~ map $\psi\colon \Dh \to A_{\infty} \cap
A'$ such that
\begin{enumerate}
\item[(i)] $b \psi(\be_{\Dh}) = b$ for all $b \in A$.
\item[(ii)] $b(\psi(dd') - \psi(d) \psi(d')) = 0 $ for all $b \in A$ and all $d,d'
\in \Dh$.
\end{enumerate}
\end{enumerate}
\end{props}

\begin{nproof}
a) $\Rightarrow$ c): By \cite[Proposition~1.9]{TomsWinter:ssa}
there is a sequence of $^{*}$-homomorphisms $\varphi_{n}\colon \Dh
\otimes \Dh \to \Dh$ such that $\varphi_{n}(d \otimes \be_{\Dh})
\to d$ for all $d \in \Dh$ as $n \to \infty$. Identify $A$ with
$A\otimes \Dh$.  Define a sequence of $^{*}$-homomorphisms
\[
\kappa_{n}\colon  A \otimes \Dh \to A \otimes \Dh
\]
by
\[
\kappa_{n}:= (\id_{A} \otimes \varphi_{n})
\circ (\id_{A} \otimes \id_{\Dh} \otimes \be_{\Dh}).
\]
Let $A^{\sim}$ denote the smallest unitization of $A$; it is clear that
$A \otimes \Dh$ is an essential ideal in $A^{\sim} \otimes \Dh$,
whence the inclusion of $A \otimes \Dh$ extends to a unital  embedding
$\iota\colon A^{\sim} \otimes \Dh \hookrightarrow \Mh(A \otimes \Dh)$.
We thus have unital $^{*}$-homomorphisms
\[
\mu_{n}\colon  \Dh \to \Mh(A \otimes \Dh)
\]
given by
\[
\mu_{n}:=  (\id_{A^{\sim}} \otimes \varphi_{n})
\circ(\be_{A^{\sim}} \otimes \be_{\Dh} \otimes \id_{\Dh}).
\]
The maps $\kappa:=\kappa_{n_{0}}$ and $\mu:=\mu_{n_{0}}$
for some large enough $n_{0}$ will have the right properties.

c) $\Rightarrow$ d): Let $b_1,b_2,\dots \in A$ and $d_1,d_2,\dots
\in \Dh$ be dense sequences in the unit balls of these algebras.
Let
 $\kappa_n \colon A
\to A$ and $\mu_n \colon \Dh \to \Mh(A)$ be sequences of
homomorphisms which satisfy
\[
 [\kappa_n(A),\mu_n(\Dh)]=0  \mbox{ and } \|\kappa_n(b_k)-b_k\| < 1/n
\]
for all $k \leq n$. Using a (possibly uncountable)  approximate
unit for $A$ which is quasicentral for $\Mh(A)$, and since
$\{\mu_{n}(d_{k}) \mid n,k \in \N\}$ is countable, it is routine
to find a countable approximate unit  $h_1,h_2,\dots$ for
$A$ such that $\|b_kh_n - b_k\|<1/n$
for all $k\leq n$, and such that $\|[\mu_n(d_k),h_n]\| < 1/n$ for
all $k \leq n$. Now, define a c.p.c.~ map $\psi \colon \Dh \to
\prod_{\N} A/\bigoplus_{\N}A$ by $d \mapsto \{h_1 \mu_1(d) h_1, h_2
\mu_2(d) h_2, \dots\}$ (i.e., the equivalence class of this
sequence), and it is easy to check that $\psi$ has the desired
properties.

d) $\Rightarrow$ b): Choose a c.p.c.~ map $\psi \colon \Dh \to
A_{\infty} \cap A'$. Use the Choi--Effros Theorem to pick some
c.p.c.~ lift $\tilde{\psi} \colon \Dh \to \prod_{\N} A$. Denote by
$\psi_n$ the $n$'th component of this map (i.e. the composition of
the projection onto the $n$'th summand of $\prod_{\N}A$ with
$\tilde{\psi}$). It is easy to verify now that for a sufficiently
large $n$, $\psi_n \colon \Dh \to A$ will satisfy the required
conditions.

b) $\Rightarrow$ a): By separability of $A$ and $\Dh$ it is
straightforward to construct a  $^{*}$-homomorphism
$\overline{\psi}\colon  A \otimes \Dh \to \prod_{\N}
A/\bigoplus_{\N}A$ such that $\overline{\psi}(a \otimes
\be_{\Dh})=a$ for all  $a \in A$. This implies $A \cong A \otimes
\Dh$ by \cite[Theorem~2.3]{TomsWinter:ssa}.
\end{nproof}
\en

\bn \label{good-maps-def} To establish $\Dh$-stability in Theorem
\ref{D-stable-bundles}, we shall mainly make use of
characterization b) in Proposition
\ref{D-stability-characterization} above. Most of the work will go
into establishing the special case $X = [0,1]$. To this end, we
introduce some ad-hoc terminology. Suppose $A$ is a separable
$\Ch([0,1])$-algebra.

\begin{dfs} Suppose we have $\varepsilon>0$, and finite sets
$\Fh \subset A$, $\Gh \subset \Dh$. Suppose $I \subset [0,1]$ is a
closed subset, and $\psi \colon \Dh \rightarrow A$ is a c.p.c.~
map. We shall say that $\psi$ is $(\Fh ;
\Gh,\varepsilon)$-\emph{good for $I$} if it satisfies the
following conditions (U),(C) and (M).
\begin{enumerate}
\item[(U)] $\left \| (b \psi(\be_{\Dh})  - b)_{I} \right \| < \varepsilon$
for all $b \in \Fh$
\item[(C)] $\|(b \psi(d) - \psi(d) b)_I\| < \varepsilon$ for all
$b \in \Fh$ and  $d \in \Gh$
\item[(M)] $\left \|(b(\psi(dd') - \psi(d) \psi(d')))_{I} \right \|< \varepsilon$
for all $b \in \Fh$ and  $d,d' \in \Gh$.
\end{enumerate}
Suppose we have some other finite set $\Gh' \subset \Dh$ and some
$\varepsilon'>0$. We shall say that $\psi$ is $(\Fh ;
\Gh,\varepsilon ; \Gh',\varepsilon')$-\emph{good for $I$}, where
$I$ here is an interval, if $\psi$ is $(\Fh ; \Gh,\varepsilon)$
good for $I$, and there is some closed neighborhood $V$ of the
endpoints of $I$ such that $\psi$ is $(\Fh ;
\Gh',\varepsilon')$-good for $V$.
\end{dfs}

The notation ``(U),(C),(M)'' is intended to serve as mnemonic for
(almost) `unital', `central' and `multiplicative', respectively.
In our case, we shall have some auxiliary $\Gh \subseteq \Gh'$ and
$\varepsilon \geq \varepsilon'$, so saying that $\psi$ is $(\Fh ;
\Gh,\varepsilon ; \Gh',\varepsilon')$-good for $I$ should be
thought of as saying that $\psi$ is $(\Fh;\Gh,\varepsilon)$-good
for $I$, and even `better' near the endpoints.
  Note that
Proposition \ref{D-stability-characterization} b) asserts that $A$
is $\Dh$-stable if and only if for any $\varepsilon>0$ and finite
sets $\Fh \subset A$, $\Gh \subset \Dh$, there is a c.p.c.~ map
$\psi:\Dh \to A$ which is $(\Fh ; \Gh , \varepsilon)$-good for
$[0,1]$.
 \en

\bn \label{commuting-good-maps}
\begin{lms}
Let $\Dh$ be a strongly self-absorbing $C^*$-algebra. Let $A$ be a
separable $\Ch([0,1])$-algebra such that $A_x$ is $\Dh$-stable for
all $x \in [0,1]$. Given $\varepsilon>0$ and finite subsets $\Fh
\subset A$ and  $\be_{\Dh} \in \Gh \subset \Dh$, there exist an
$n\in \N$, c.p.c.~ maps $\psi_1,\dots,\psi_n \colon \Dh \to A$,
and points $0=t_0<t_1<\dots<t_n = 1$ such that $\psi_k$ is
$(\Fh;\Gh,\varepsilon)$-good for $[t_{k-1},t_k]$ for $k=1,...,n$.
\end{lms}

\begin{nproof}
For any $x \in [0,1]$ we can find, by Proposition
\ref{D-stability-characterization} b), a c.p.c.~ map $\sigma_x
\colon \Dh \to A_x$ such that for all $b \in \Fh$, $d,d' \in \Gh$,
we have
\begin{enumerate}
\item[(U)] $\|b_x \sigma_x(\be_{\Dh}) - b_x\| < \varepsilon$
\item[(C)] $\|b_x \sigma_x(d) - \sigma_x(d) b_x\| < \varepsilon $
\item[(M)] $\| b_x(\sigma_x(dd') - \sigma_x(d) \sigma_x(d'))\|< \varepsilon$.
\end{enumerate}
Use the Choi--Effros theorem to find a c.p.c.~ lift $\rho_x \colon
\Dh \to A$ of $\sigma_x$. By upper semicontinuity of the norm
function (see the discussion in \ref{upper-semicontinuity} above),
it follows that for any $x$ there is some open interval $I_x$
containing $x$ such that $\rho_x$ is $(\Fh ; \Gh,
\varepsilon)$-good for $\overline{I_x}$. Those intervals cover
$[0,1]$, so by compactness, we can find a finite subcover. Now, by
making the intervals smaller, and by omitting redundant elements
of the cover, we obtain $\psi_1,...,\psi_n$ as required.
\end{nproof}
\en

We shall prove Theorem \ref{D-stable-bundles} for $X=[0,1]$ by
patching the maps $\psi_1,..,\psi_n$ from Lemma
\ref{commuting-good-maps} (which will be chosen with a finer
auxiliary $\Gh',\varepsilon'$) into one c.p.c.~ map. In the
following two lemmas, we will show how two c.p.c.~ maps
$\rho,\sigma \colon \Dh \to A$, which are sufficiently well-behaved on adjacent
intervals $[r,s]$ and $[s,t]$ respectively, can be patched
together to obtain one c.p.c.~ map which is well-behaved on the
union (where `well-behaved' here means $(\Fh ; \Gh,\varepsilon ;
\Gh' , \varepsilon')$-good for some given
$\Fh,\Gh,\varepsilon,\Gh',\varepsilon'$). Since those lemmas are
somewhat technical, we first give here a hand-waved outline of the
proof.

In Lemma \ref{mu-lemma}, we will show that, after perturbing
$\rho$ and $\sigma$, we can find  auxiliary c.p.c.~ maps
$\nu_{\rho},\nu_{\sigma} \colon \Dh \to A$ which are well-behaved
near the point $s$, and such that near $s$, we have that
$\nu_{\rho} \approx \nu_{\sigma}$, that the range of $\nu_{\rho}$
approximately commutes with the range of $\rho$, and that the
range of $\nu_{\sigma}$ approximately commutes with the range of
$\sigma$ -- in fact, we obtain c.p.c.~ maps
$\mu_{\rho},\mu_{\sigma} \colon \Dh \otimes \Dh \to A$ such that
$\mu_{\rho} \approx \mbox{``}\rho \otimes \nu_{\rho} \mbox{''}$
and $\mu_{\sigma} \approx \mbox{``}\sigma \otimes \nu_{\sigma}
\mbox{''}$ (since the maps do not exactly commute, we cannot
actually define $\mbox{``}\rho \otimes \nu_{\rho} \mbox{''} ,
\mbox{``}\sigma \otimes \nu_{\sigma} \mbox{''}$ as c.p.c.~ maps).

 In Lemma \ref{patching
lemma}, we shall construct a new c.p.c.~ map which is well-behaved
on $[r,t]$. To do so, we take a unitary path $\{u_x\}_{x \in
[0,1]}\subset \Dh \otimes \Dh$, such that $u_0 = \be_{\Dh \otimes
\Dh}$ and $u_1 (d \otimes \be_{\Dh}) u_1^* \approx \be_{\Dh}
\otimes d$. We then have that $\mu_{\rho}(u_0 (d
\otimes\be_{\Dh})u_0^*) \approx \rho(d)$, and $\mu_{\rho}(u_1 (d
\otimes\be_{\Dh})u_1^*) \approx \nu_{\rho}(d) \approx
\nu_{\sigma}(d)$. So, in this way, we can `connect' $\rho$ to
$\nu_{\rho} (\approx \nu_{\sigma})$ `along' the path of unitaries,
and then similarly `connect' $\nu_{\sigma}$ to $\sigma$. We thus
obtain the desired c.p.c.~ map $\psi$, which agrees with $\rho$
near $r$, with $\sigma$ near $t$, and is roughly $\nu_{\rho}$
($\approx \nu_{\sigma}$) near the midpoint $s$.

We finally note that if we were to restrict ourselves to unital
$C^*$-algebras, the computations would become somewhat simpler.
The reader who is happy to make this assumption in first reading
should replace ``c.p.c.'' by ``u.c.p.'' throughout (one can then
omit condition (U), and omit the ``b'' from condition (M)).

\bn \label{mu-lemma}
\begin{lms}
Let $\Dh$ be a strongly self-absorbing $C^*$-algebra, and let $A$
be a separable $\Ch([0,1])$-algebra. Suppose we have
$\varepsilon>0$ and finite sets $\Fh \subset A$, $\be_{\Dh} \in
\Gh \subset \Dh$ consisting of elements of norm at most one, with
$\Fh = \Fh^*$, $\Gh = \Gh^*$. Suppose we have points $0 \leq r < s
< t \leq 1$ and two c.p.c.~ maps $\rho,\sigma \colon \Dh \to A$
which are $(\Fh ; \Gh , \varepsilon)$-good for $[r,s]$, $[s,t]$
respectively. Suppose furthermore that $A_s$ is $\Dh$-stable.

It follows that there are c.p.c.~ maps $\rho' , \sigma' \colon \Dh
\to A$ which are $(\Fh ; \Gh , \varepsilon)$-good for $[r,s]$,
$[s,t]$ respectively, and c.p.c.~ maps $\nu_{\rho'},\nu_{\sigma'}
\colon \Dh \to A$, $\mu_{\rho'},\mu_{\sigma'} \colon \Dh \otimes
\Dh \to A$ such that  $\nu_{\rho'}$ and $\nu_{\sigma'}$ are $(\Fh
; \Gh , 3\varepsilon)$-good for some interval $I\subset (r,t)$
containing $s$ in its interior, and such that for any $b \in \Fh$
and any $d,d' \in \Gh$ we have
\begin{enumerate}
\item[(i)$_{\rho}$] $\|(b[\rho'(d),\nu_{\rho'}(d')])_{I}\| <
2\varepsilon$
\item[(i)$_{\sigma}$]
$\|(b[\sigma'(d),\nu_{\sigma'}(d')])_{I}\| < 2\varepsilon$
\item[(ii)$_{\rho}$] $\|(b(\rho'(d)\nu_{\rho'}(d') - \mu_{\rho'}(d \otimes d'))_{I}\| <
\varepsilon$
\item[(ii)$_{\sigma}$] $\|(b(\sigma'(d)\nu_{\sigma'}(d') - \mu_{\sigma'}(d \otimes d'))_{I}\| <
\varepsilon$
\item[(iii)$_{\; }$] $\|(b(\nu_{\rho'}(d) -
\nu_{\sigma'}(d)))_{I}\| < 2\varepsilon$.
\end{enumerate}
If $\rho,\sigma$ are $(\Fh ; \Gh , \varepsilon ; \Gh' ,
\varepsilon')$-good for $[r,s]$,$[s,t]$ respectively, for some
finite self-adjoint set $\Gh' \supset \Gh$ in the unit ball of
$\Dh$ and for some $0 < \varepsilon' < \varepsilon$, then one can
arrange that so are $\rho',\sigma'$, that $\nu_{\rho}$ and
$\nu_{\sigma}$ are $(\Fh ; \Gh' , 3\varepsilon')$-good for the
interval $I$, and that conditions
(i)$_{\rho}$,(i)$_{\sigma}$,(ii)$_{\rho}$,(ii)$_{\sigma}$,(iii)
hold with $\varepsilon'$ instead of $\varepsilon$, and $\Gh'$
instead of $\Gh$.
\end{lms}

\begin{nproof}
Denote $\Gh\cdot\Gh = \{ab \mid a,b \in \Gh\}$, then
$\Gh\cdot\Gh \supset \Gh$. Since $\Fh,\Gh$ are finite, it follows
that there is some $\varepsilon > \eta > 0$ such that the maps
$\rho,\sigma$ are also $(\Fh;\Gh,\eta)$-good for $[r,s]$,$[s,t]$,
respectively. We use Proposition
\ref{D-stability-characterization} c) to find $^*$-homomorphisms
$\kappa \colon A_s \to A_s$ and $\mu:\Dh \to \Mh(A_s)$ such that
\begin{equation}
\label{kappa-1}
 [\kappa(A_s),\mu(\Dh)]=0
 \end{equation}
 and
 \begin{equation}
 \label{kappa-2}
 \|\kappa(a_{s})-a_{s}\| < (\varepsilon-\eta)/3
\end{equation}
for all $a \in \Fh \cup \rho(\Gh\cdot\Gh) \cup
\sigma(\Gh\cdot\Gh)$. This in particular implies that
\begin{equation} \label{eq-mu-comm}
 \|[\mu(d),b_s]\| < \|[\mu(d),\kappa(b_s)]\| + 2 (\varepsilon-\eta)/3 =
 2 (\varepsilon-\eta)/3 <
\varepsilon
\end{equation}
 for all $b \in \Fh$, $d \in \Gh$. Use the Choi--Effros
theorem to find c.p.c.~ lifts $\tilde{\rho},\tilde{\sigma} \colon
\Dh \to A$ for the maps $\kappa \circ \pi_{s} \circ \rho$,
$\kappa \circ \pi_{s} \circ \sigma$, respectively. Using (\ref{kappa-2}),
we get

$ \begin{array}{ll} \vspace{5pt}
 \mbox{(U)}_{\kappa \circ
\sigma} & \displaystyle  \|b_s \kappa(\sigma(\be_{\Dh})_s) - b_s\|
< \eta + \frac{\varepsilon-\eta}{3}  < \varepsilon
 \\ \vspace{5pt}
\mbox{(C)}_{\kappa \circ \sigma} & \displaystyle \|b_s \kappa
(\sigma(d)_s) - \kappa (\sigma(d)_s) b_s\|
 \leq \|[b_s,\sigma(d)_s]\|
  +  2\cdot \|\sigma(d)_s -
 \kappa(\sigma(d)_s)\|
  \\  &
  \displaystyle
  <  \eta + 2 \cdot \frac{\varepsilon-\eta}{3}  < \varepsilon
  \\
\mbox{(M)}_{\kappa \circ \sigma} & \displaystyle \|
b_s(\kappa(\sigma(dd')_s) - \kappa(\sigma(d)_s)\kappa(
\sigma(d')_s))\|
   < \eta + 3 \cdot
\frac{\varepsilon-\eta}{3}  = \varepsilon
\end{array} $

\noindent for all $b \in \Fh$, $d,d' \in \Gh$ and the same estimates hold
for $\rho$ instead of
$\sigma$. Thus, by upper semicontinuity  of the norm function, we
have that there is some $\gamma>0$ such that
$\tilde{\rho},\tilde{\sigma}$ are $(\Fh; \Gh,\varepsilon)$-good
for $I:= [s-\gamma,s+\gamma]$, and such that
\begin{equation} \label{rho-rhotilde}
 \|(\rho(d) - \tilde{\rho}(d))_{[s-\gamma,s+\gamma]}\| \,,\,
 \|(\sigma(d) - \tilde{\sigma}(d))_{[s-\gamma,s+\gamma]}\| <
 \frac{\varepsilon-\eta}{3}
\end{equation}
 for all $d\in \Gh\cdot\Gh$. We may assume
that $r < s-\gamma < s +\gamma < t$. Define $f,g \in \Ch([0,1])$
as follows.
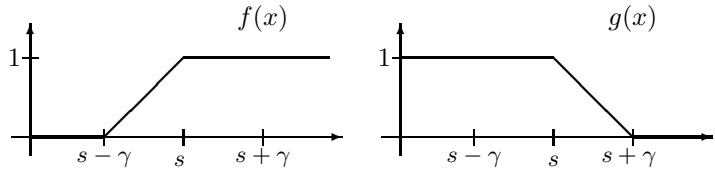
\begin{center}
\begin{picture}(290,65)
 \put(0,10){\vector(1,0){125}}
 \put(7,3){\vector(0,1){50}}
 \put(65,7){\line(0,1){6}}
 \put(35,7){\line(0,1){6}}
 \put(95,7){\line(0,1){6}}
 \put(4,40){\line(1,0){6}}
\thicklines
  \put(7,10){\line(1,0){28}}
 \put(35,10){\line(1,1){30}}
 \put(65,40){\line(1,0){55}}
 \put(1,40){\makebox(0,0){$1$}}
 \put(65,-1){\makebox(0,0)[b]{\small $s$ \normalsize}}
 \put(35,-1){\makebox(0,0)[b]{\small $s-\gamma$\normalsize}}
 \put(95,-1){\makebox(0,0)[b]{\small $s+\gamma$\normalsize}}
 \put(95,55){\makebox(0,0){$f(x)$}}

 \thinlines

 \put(140,10){\vector(1,0){125}}
 \put(147,3){\vector(0,1){50}}
 \put(235,7){\line(0,1){6}}
 \put(205,7){\line(0,1){6}}
 \put(175,7){\line(0,1){6}}
 \put(144,40){\line(1,0){6}}
\thicklines
 \put(147,40){\line(1,0){58}}
 \put(205,40){\line(1,-1){30}}
 \put(235,10){\line(1,0){30}}
 \put(141,40){\makebox(0,0){$1$}}
 \put(235,-1){\makebox(0,0)[b]{\small $s+\gamma$ \normalsize}}
 \put(205,-1){\makebox(0,0)[b]{\small $s$\normalsize}}
 \put(175,-1){\makebox(0,0)[b]{\small $s-\gamma$\normalsize}}
 \put(235,55){\makebox(0,0){$g(x)$}}
\end{picture}
\end{center}

Now, define
\[
\rho'(d) := (1-f) \cdot \rho(d) + f \cdot
\tilde{\rho}(d) \mbox{ and } \sigma'(d) := (1-g) \cdot \sigma(d) + g \cdot
\tilde{\sigma}(d).
\]
It is clear that $\rho',\sigma'$ are c.p.c.
 We now show that $\rho'$ is
$(\Fh;\Gh,\varepsilon)$-good for $[r,s]$ and that $\sigma'$ is
$(\Fh;\Gh,\varepsilon)$-good for $[s,t]$. We establish this for
$\rho'$ -- the proof for $\sigma'$ is the same.

Note that $\rho'(d)_{[0,s-\gamma]} =  \rho(d)_{[0,s-\gamma]}$, and
thus it remains to show that $\rho'$ is
$(\Fh;\Gh,\varepsilon)$-good for $[s-\gamma,s]$. Indeed, using
(\ref{rho-rhotilde}), we obtain

$
\begin{array}{ll} \vspace{5pt}
\mbox{(U)}_{\rho'} &  \displaystyle \|(b \rho'(\be_{\Dh}) -
b)_{[s-\gamma,s]}\|
 \\\vspace{5pt} &  \displaystyle < \|(b \rho(\be_{\Dh}) - b)_{[s-\gamma,s]}\| +
\frac{\varepsilon-\eta}{3} < \eta + \frac{\varepsilon-\eta}{3} <
\varepsilon
 \\ \vspace{5pt}
\mbox{(C)}_{\rho'} & \displaystyle  \|(b \rho'(d) - \rho'(d)
b)_{[s-\gamma,s]}\|
 \\ \vspace{5pt}
 & < \displaystyle
 \|(b \rho(d) - \rho(d) b)_{[s-\gamma,s]}\| + 2 \cdot\frac{\varepsilon-\eta}{3}
 <\eta + 2\cdot \frac{\varepsilon-\eta}{3} < \varepsilon
 \\ \vspace{5pt}
\mbox{(M)}_{\rho'} & \displaystyle \| (b(\rho'(dd')
 -\rho'(d)\rho'(d')))_{[s-\gamma,s]}\|
 \\ & \displaystyle < \| (b(\rho(dd')
-\rho(d)\rho(d')))_{[s-\gamma,s]}\| + 3 \cdot
\frac{\varepsilon-\eta}{3} <
 \eta + \varepsilon-\eta = \varepsilon
\end{array} $

\noindent for all $b \in \Fh$, $d,d' \in \Gh$, as required.

Now, note that by (\ref{kappa-1}), the range of $\mu$ commutes
with $\rho'(d)_s$ and $\sigma'(d)_s$ for all $d \in \Dh$. We may
thus define c.p.c.~ maps $\bar{\mu}_{\rho'},\bar{\mu}_{\sigma'}
\colon \Dh \otimes \Dh \to A_s$ by
$$\bar{\mu}_{\rho'}(d \otimes d') :=
\rho'(d)_s\mu(d') \;\mbox{ and }\; \bar{\mu}_{\sigma'}(d \otimes
d') := \sigma'(d)_s\mu(d')
$$
and we define c.p.c.~ maps $\bar{\nu}_{\rho'},\bar{\nu}_{\sigma'}
\colon \Dh \to A_s$ by
$$\bar{\nu}_{\rho'}(d) :=
\rho'(\be_{\Dh})_s\mu(d') = \bar{\mu}_{\rho'}(\be_{\Dh} \otimes
d') \,\mbox{ and }\, \bar{\nu}_{\sigma'}(d) :=
\sigma'(\be_{\Dh})_s\mu(d') = \bar{\mu}_{\sigma'}(\be_{\Dh}
\otimes d').
$$
We now use the Choi--Effros theorem to choose some c.p.c.~ lifts
$\mu_{\rho'},\mu_{\sigma'},\nu_{\rho'},\nu_{\sigma'}$ into $A$ for
$\bar{\mu}_{\rho'},\bar{\mu}_{\sigma'},\bar{\nu}_{\rho'},\bar{\nu}_{\sigma'}$.
We need to show that $\nu_{\rho'}$ and $\nu_{\sigma'}$ are $(\Fh ;
\Gh , 3\varepsilon)$-good for some neighborhood of $s$, and that
the five conditions in the statement hold. By upper
semicontinuity, it suffices to verify this at $s$.

We first show that $\nu_{\rho'}$ is $(\Fh ; \Gh ,
3\varepsilon)$-good for some neighborhood of $s$; the proof for
$\nu_{\sigma'}$ is the same, so we omit it. Indeed, since $\rho'$
is $(\Fh;\Gh,\varepsilon)$-good for $\{s\}$, we have

$
\begin{array}{ll} \vspace{5pt}
\mbox{(U)}_{\nu_{\rho'}} & \displaystyle
\|b_s\nu_{\rho'}(\be_{\Dh})_s-b_s\| =
 \|b_s\rho'(\be_{\Dh})_s -
b_s\| < \varepsilon
 \\\vspace{5pt}
 \mbox{(C)}_{\nu_{\rho'}} &
 \displaystyle \|b_s \nu_{\rho'}(d)_s - \nu_{\rho'}(d)_sb_s\|
  = \|b_s \rho'(\be_{\Dh})_s\mu(d) -
\underbrace{\rho'(\be_{\Dh})_s \mu(d)}_{=\mu(d)\rho'(\be_{\Dh})_s}
b_s\|
 \\\vspace{5pt} &
 \leq
 \underbrace{\|b_s\rho'(\be_{\Dh})_s-b_s\|}_{< \varepsilon \;
 \textrm{  by (U)}_{\rho'}}
 + \underbrace{\|\rho'(\be_{\Dh})_sb_s-b_s\|}_{< \varepsilon \;
 \textrm{  by (U)}_{\rho'}}
 +  \underbrace{\|b_s \mu(d) - \mu(d)b_s\|}_{< \varepsilon \;
 \textrm{  by (\ref{eq-mu-comm})}} < 3\varepsilon
 \\\vspace{5pt}
\mbox{(M)}_{\nu_{\rho'}} & \displaystyle \|b_s (\nu_{\rho'}(dd')_s
- \nu_{\rho'}(d)_s
 \nu_{\rho'}(d')_s)\|
 \\\vspace{5pt} &  =
  \|b_s\rho'(\be_{\Dh})_{s}(\mu(dd') -
  \underbrace{\mu(d)\rho'(\be_{\Dh})_{s}}_{\mbox{}=
  \rho'(\be_{\Dh})_{s}\mu(d)}
  \mu(d')\|
  \\\vspace{5pt}  &  =
  \|b_s \rho'(\be_{\Dh})_s(\underbrace{\mu(dd') -
  \mu(d)\mu(d')}_{\mbox{}= 0})
  \\\vspace{5pt} & \hspace{1em}  +
  (b_s
  -b_s\rho'(\be_{\Dh})_s)\rho'(\be_{\Dh})_s\mu(d)\mu(d')\|
  \\\vspace{5pt}  &  \leq \|b_s-b_s\rho'(\be_{\Dh})_s\|<
  \varepsilon.
\end{array} $

\noindent We now verify condition (i)$_{\rho}$:
\begin{eqnarray*}
\lefteqn{ \|(b_s[\rho'(d)_s,\nu_{\rho'}(d')_s]) \|
  =
 \| b_s[\rho'(d)_s,\rho'(\be_{\Dh})_s]\mu(d')\|}
   \\  &\mbox{}\leq&
   \underbrace{\|b_s(\rho'(d)_s\rho'(\be_{\Dh})_s -
   \rho'(d)_s)\|}_{< \varepsilon \;\textrm{ by
   (M)}_{\rho'}}
 + \underbrace{\|b_s(\rho'(\be_{\Dh})_s\rho'(d)_s -
 \rho'(d)_s)\|}_{< \varepsilon \;\textrm{ by
   (M)}_{\rho'}}
< 2\varepsilon \,.
\end{eqnarray*}
As for condition (ii)$_{\rho}$, we have indeed (again using
(M)$_{\rho'}$):
\begin{eqnarray*}
\|(b_s(\rho'(d)_s\nu_{\rho'}(d')_s - \mu_{\rho'}(d \otimes
d')_s)\|
 = \|(b_s(\rho'(d)_s\rho'(\be_{\Dh})_s - \rho'(d)_s) \mu(d')\| < \varepsilon.
\end{eqnarray*}
 Conditions
(i)$_{\sigma}$ and (ii)$_{\sigma}$ follow in the same manner. It
remains to check condition (iii). Indeed, using condition (U) for
$\rho'$ and $\sigma'$, we have
\begin{eqnarray*}
 \lefteqn{\|(b_s(\nu_{\rho'}(d)_s - \nu_{\sigma'}(d)_s)\|
 = \|b_s(\rho'(\be_{\Dh})_s - \sigma'(\be_{\Dh})_s)\mu(d)\| }
 \\ & \leq &
 \|b_s\rho'(\be_{\Dh})_s - b_s\| + \|b_s - b_s\sigma'(\be_{\Dh})_s\|
 < 2\varepsilon \,.
\end{eqnarray*}

The last part of the lemma follows immediately from the
construction, provided that we choose $\gamma>0$ such that $\rho$
and $\sigma$ are $(\Fh ; \Gh' , \varepsilon')$-good for
$[s-\gamma,s+\gamma]$.
\end{nproof}
\en

\bn \label{patching lemma}
\begin{lms}
Let $\Dh$ be a $K_{1}$-injective strongly self-absorbing
$C^{*}$-algebra. Let $A$ be a separable $\Ch([0,1])$-algebra.
Suppose we are given $\varepsilon>0$ and finite subsets $\Fh
\subset A$ and $\be_{\Dh} \in \Gh \subset \Dh$, with $\Fh =
\Fh^*$, $\Gh = \Gh^*$.

There exist $0<\varepsilon'<\varepsilon$ and a finite set $\Gh
\subset \Gh' \subset \Dh$, such that if we have two c.p.c.~ maps
$\rho,\sigma \colon \Dh \rightarrow A$, and points $0 \leq r  < s
< t \leq 1$ such that $\rho$ is $(\Fh ; \Gh,\varepsilon ;
\Gh',\varepsilon')$-good for $[r,s]$ and $\sigma$ is $(\Fh ;
\Gh,\varepsilon ; \Gh',\varepsilon')$-good for $[s,t]$,
 then
there is a c.p.c.~ map $\psi \colon \Dh  \rightarrow A$
 which is $(\Fh ; \Gh,\varepsilon ;
\Gh',\varepsilon')$-good for $[r,t]$.
\end{lms}

\begin{nproof}
We may assume without loss of generality that the elements of
$\Fh$ and $\Gh$ have norm at most one. Since $\Dh$ is
$K_1$-injective, \cite[Proposition~1.13]{TomsWinter:ssa}
guarantees that we may choose a unitary $u \in \Ch([0,1],\Dh
\otimes \Dh)$ such that
\begin{equation}
 u_{0}= \be_{\Dh \otimes \Dh}
\end{equation}
and
\begin{equation}
\label{u_1}
 \|u_{1} (d \otimes \be_{\Dh}) u_{1}^{*} - \be_{\Dh} \otimes
d\| < \frac{\varepsilon}{4}
\end{equation}
for all $d \in \Gh \cdot \Gh = \{d_{1} d_{2}\mid d_{1},d_{2} \in
\Gh\}$. Since the set $\{u_x\}_{x \in [0,1]}$ is compact, we may
now fix some $m \in \N$ and $v_{ji},w_{ji}$, $i,j=1,\dots,m$ in
the unit ball of $\Dh$, such that for any $x \in [0,1]$ there is
an $i$ with
$$
 \Big\|u_{x} - \sum_{j=1}^{m} v_{ji} \otimes w_{ji} \Big\|
< \frac{\varepsilon}{9}.
$$
Denote $y_i := \sum_{j=1}^{m} v_{ji} \otimes w_{ji} $. We may
assume that $\|y_i \| \leq 1$ for all $i$.

Now, set
$$
\Gh' := \left \{d_1d_2d_3d_4d_5d_6 \mid d_1,...,d_6 \in \Gh \cup
\{ v_{ji} , v_{ji}^* , w_{ji} , w_{ji}^* \mid i,j = 1,\dots,m\}
\right \}
$$
and
$$
\varepsilon' := \frac{\varepsilon}{144m^4} \,.
$$

Suppose we are given $\rho,\sigma$ as in the statement. We wish to
construct the desired c.p.c.~ map $\psi$. By replacing
$\rho,\sigma$ by $\rho',\sigma'$ as in the Lemma \ref{mu-lemma} if
needed, we have c.p.c.~ maps $\nu_{\rho},\nu_{\sigma} \colon \Dh
\to A$ and $\mu_{\rho},\mu_{\sigma} \colon \Dh \otimes \Dh \to A$
such that the conclusions of the lemma hold for some interval $I
\subset (r,t)$ containing $s$ in its interior. By upper
semicontinuity of the norm function, there is some $\delta>0$ such
that both $\rho$ and $\sigma$ are $(\Fh;\Gh',\varepsilon')$-good
for $[s-3\delta,s+3\delta] \subset I$.

 Define c.p.c.~ maps
$$\varphi_{\rho},\varphi_{\sigma} \colon \Ch([0,1]) \otimes \Dh \otimes \Dh \rightarrow
A$$ by
$$
\varphi_{\rho}(f \otimes d \otimes d') := f \cdot \mu_{\rho}(d
\otimes d')
 \; , \;
 \varphi_{\sigma}(f \otimes d \otimes d') := f \cdot
\mu_{\sigma}(d \otimes d') .
$$
 Define continuous functions $h_1,h_2,h_3,h_4 \colon [0,1] \to [0,1]$
 which sum up to $1$, by:
 
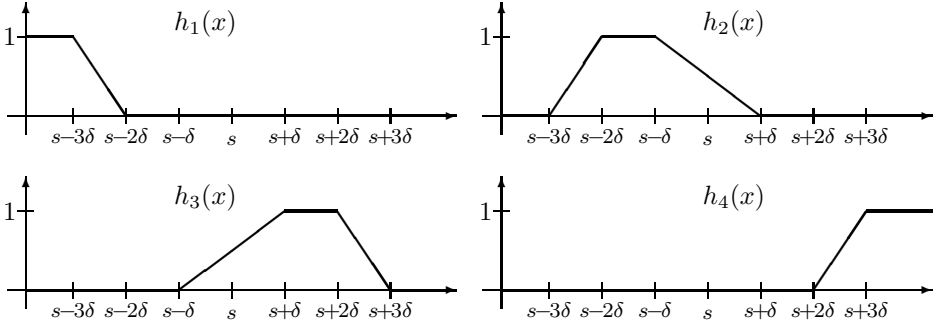
\begin{center}
\begin{picture}(400,65)
 \put(0,10){\vector(1,0){170}}
 \put(7,3){\vector(0,1){50}}
 \put(25,7){\line(0,1){6}}
 \put(45,7){\line(0,1){6}}
 \put(65,7){\line(0,1){6}}
 \put(85,7){\line(0,1){6}}
 \put(105,7){\line(0,1){6}}
 \put(125,7){\line(0,1){6}}
 \put(145,7){\line(0,1){6}}
 \put(4,40){\line(1,0){6}}
\thicklines
 \put(7,40){\line(1,0){18}}
 \put(25,40){\line(2,-3){20}}
 \put(45,10){\line(1,0){125}}
 \put(1,40){\makebox(0,0){$1$}}
 \put(25,-1){\makebox(0,0)[b]{\footnotesize $s\!\!-\!\!3\delta$\normalsize}}
 \put(45,-1){\makebox(0,0)[b]{\footnotesize  $s\!\!-\!\!2\delta$\normalsize}}
 \put(65,-1){\makebox(0,0)[b]{\footnotesize  $s\!\!-\!\!\delta$\normalsize}}
 \put(85,-1){\makebox(0,0)[b]{\footnotesize  $s$\normalsize}}
 \put(105,-1){\makebox(0,0)[b]{\footnotesize  $s\!\!+\!\!\delta$\normalsize}}
 \put(125,-1){\makebox(0,0)[b]{\footnotesize  $s\!\!+\!\!2\delta$\normalsize}}
 \put(145,-1){\makebox(0,0)[b]{\footnotesize  $s\!\!+\!\!3\delta$\normalsize}}
 \put(75,45){\makebox(0,0){$h_1(x)$}}

 \thinlines

  \put(180,10){\vector(1,0){170}}
 \put(187,3){\vector(0,1){50}}
 \put(205,7){\line(0,1){6}}
 \put(225,7){\line(0,1){6}}
 \put(245,7){\line(0,1){6}}
 \put(265,7){\line(0,1){6}}
 \put(285,7){\line(0,1){6}}
 \put(305,7){\line(0,1){6}}
 \put(325,7){\line(0,1){6}}
 \put(184,40){\line(1,0){6}}
\thicklines
 \put(187,10){\line(1,0){18}}
 \put(205,10){\line(2,3){20}}
 \put(225,40){\line(1,0){20}}
 \put(245,40){\line(4,-3){40}}
 \put(285,10){\line(1,0){65}}
 \put(181,40){\makebox(0,0){$1$}}
 \put(205,-1){\makebox(0,0)[b]{\footnotesize $s\!\!-\!\!3\delta$\normalsize}}
 \put(225,-1){\makebox(0,0)[b]{\footnotesize  $s\!\!-\!\!2\delta$\normalsize}}
 \put(245,-1){\makebox(0,0)[b]{\footnotesize  $s\!\!-\!\!\delta$\normalsize}}
 \put(265,-1){\makebox(0,0)[b]{\footnotesize  $s$\normalsize}}
 \put(285,-1){\makebox(0,0)[b]{\footnotesize  $s\!\!+\!\!\delta$\normalsize}}
 \put(305,-1){\makebox(0,0)[b]{\footnotesize  $s\!\!+\!\!2\delta$\normalsize}}
 \put(325,-1){\makebox(0,0)[b]{\footnotesize  $s\!\!+\!\!3\delta$\normalsize}}
 \put(275,45){\makebox(0,0){$h_2(x)$}}

\end{picture}
\begin{picture}(400,65)
 \put(0,10){\vector(1,0){170}}
 \put(7,3){\vector(0,1){50}}
 \put(25,7){\line(0,1){6}}
 \put(45,7){\line(0,1){6}}
 \put(65,7){\line(0,1){6}}
 \put(85,7){\line(0,1){6}}
 \put(105,7){\line(0,1){6}}
 \put(125,7){\line(0,1){6}}
 \put(145,7){\line(0,1){6}}
 \put(4,40){\line(1,0){6}}
\thicklines
 \put(7,10){\line(1,0){58}}
 \put(65,10){\line(4,3){40}}
 \put(105,40){\line(1,0){20}}
 \put(125,40){\line(2,-3){20}}
 \put(145,10){\line(1,0){25}}
 \put(1,40){\makebox(0,0){$1$}}
 \put(25,-1){\makebox(0,0)[b]{\footnotesize $s\!\!-\!\!3\delta$\normalsize}}
 \put(45,-1){\makebox(0,0)[b]{\footnotesize  $s\!\!-\!\!2\delta$\normalsize}}
 \put(65,-1){\makebox(0,0)[b]{\footnotesize  $s\!\!-\!\!\delta$\normalsize}}
 \put(85,-1){\makebox(0,0)[b]{\footnotesize  $s$\normalsize}}
 \put(105,-1){\makebox(0,0)[b]{\footnotesize  $s\!\!+\!\!\delta$\normalsize}}
 \put(125,-1){\makebox(0,0)[b]{\footnotesize  $s\!\!+\!\!2\delta$\normalsize}}
 \put(145,-1){\makebox(0,0)[b]{\footnotesize  $s\!\!+\!\!3\delta$\normalsize}}
 \put(75,45){\makebox(0,0){$h_3(x)$}}

 \thinlines

 \put(180,10){\vector(1,0){170}}
 \put(187,3){\vector(0,1){50}}
 \put(205,7){\line(0,1){6}}
 \put(225,7){\line(0,1){6}}
 \put(245,7){\line(0,1){6}}
 \put(265,7){\line(0,1){6}}
 \put(285,7){\line(0,1){6}}
 \put(305,7){\line(0,1){6}}
 \put(325,7){\line(0,1){6}}
 \put(184,40){\line(1,0){6}}
\thicklines
 \put(187,10){\line(1,0){118}}
 \put(305,10){\line(2,3){20}}
 \put(325,40){\line(1,0){25}}
 \put(181,40){\makebox(0,0){$1$}}
 \put(205,-1){\makebox(0,0)[b]{\footnotesize $s\!\!-\!\!3\delta$\normalsize}}
 \put(225,-1){\makebox(0,0)[b]{\footnotesize  $s\!\!-\!\!2\delta$\normalsize}}
 \put(245,-1){\makebox(0,0)[b]{\footnotesize  $s\!\!-\!\!\delta$\normalsize}}
 \put(265,-1){\makebox(0,0)[b]{\footnotesize  $s$\normalsize}}
 \put(285,-1){\makebox(0,0)[b]{\footnotesize  $s\!\!+\!\!\delta$\normalsize}}
 \put(305,-1){\makebox(0,0)[b]{\footnotesize  $s\!\!+\!\!2\delta$\normalsize}}
 \put(325,-1){\makebox(0,0)[b]{\footnotesize  $s\!\!+\!\!3\delta$\normalsize}}
 \put(275,45){\makebox(0,0){$h_4(x)$}}

\end{picture}
\end{center}
Define $g_{\rho},g_{\sigma} \colon [0,1] \to [0,1]$ by
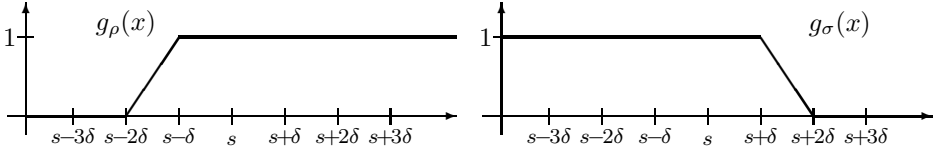
\begin{center}
\begin{picture}(380,65)
 \put(0,10){\vector(1,0){170}}
 \put(7,3){\vector(0,1){50}}
 \put(25,7){\line(0,1){6}}
 \put(45,7){\line(0,1){6}}
 \put(65,7){\line(0,1){6}}
 \put(85,7){\line(0,1){6}}
 \put(105,7){\line(0,1){6}}
 \put(125,7){\line(0,1){6}}
 \put(145,7){\line(0,1){6}}
 \put(4,40){\line(1,0){6}}
\thicklines
 \put(7,10){\line(1,0){38}}
 \put(45,10){\line(2,3){20}}
 \put(65,40){\line(1,0){105}}
 \put(1,40){\makebox(0,0){$1$}}
 \put(25,-1){\makebox(0,0)[b]{\footnotesize $s\!\!-\!\!3\delta$\normalsize}}
 \put(45,-1){\makebox(0,0)[b]{\footnotesize  $s\!\!-\!\!2\delta$\normalsize}}
 \put(65,-1){\makebox(0,0)[b]{\footnotesize  $s\!\!-\!\!\delta$\normalsize}}
 \put(85,-1){\makebox(0,0)[b]{\footnotesize  $s$\normalsize}}
 \put(105,-1){\makebox(0,0)[b]{\footnotesize  $s\!\!+\!\!\delta$\normalsize}}
 \put(125,-1){\makebox(0,0)[b]{\footnotesize  $s\!\!+\!\!2\delta$\normalsize}}
 \put(145,-1){\makebox(0,0)[b]{\footnotesize  $s\!\!+\!\!3\delta$\normalsize}}
 \put(45,45){\makebox(0,0){$g_{\rho}(x)$}}

 \thinlines

 \put(180,10){\vector(1,0){170}}
 \put(187,3){\vector(0,1){50}}
 \put(205,7){\line(0,1){6}}
 \put(225,7){\line(0,1){6}}
 \put(245,7){\line(0,1){6}}
 \put(265,7){\line(0,1){6}}
 \put(285,7){\line(0,1){6}}
 \put(305,7){\line(0,1){6}}
 \put(325,7){\line(0,1){6}}
 \put(184,40){\line(1,0){6}}
\thicklines
 \put(187,40){\line(1,0){98}}
 \put(285,40){\line(2,-3){20}}
 \put(305,10){\line(1,0){45}}
 \put(181,40){\makebox(0,0){$1$}}
 \put(205,-1){\makebox(0,0)[b]{\footnotesize $s\!\!-\!\!3\delta$\normalsize}}
 \put(225,-1){\makebox(0,0)[b]{\footnotesize  $s\!\!-\!\!2\delta$\normalsize}}
 \put(245,-1){\makebox(0,0)[b]{\footnotesize  $s\!\!-\!\!\delta$\normalsize}}
 \put(265,-1){\makebox(0,0)[b]{\footnotesize  $s$\normalsize}}
 \put(285,-1){\makebox(0,0)[b]{\footnotesize  $s\!\!+\!\!\delta$\normalsize}}
 \put(305,-1){\makebox(0,0)[b]{\footnotesize  $s\!\!+\!\!2\delta$\normalsize}}
 \put(325,-1){\makebox(0,0)[b]{\footnotesize  $s\!\!+\!\!3\delta$\normalsize}}
 \put(315,45){\makebox(0,0){$g_{\sigma}(x)$}}

\end{picture}
\end{center}

Define  unitaries $u_{\rho},u_{\sigma} \in \Ch([0,1]) \otimes \Dh
\otimes \Dh \cong \Ch([0,1],\Dh \otimes \Dh)$ by
$$
 u_{\rho\,x} := u_{g_{\rho}(x)} \; , \; u_{\sigma\,x} :=  u_{g_{\sigma}(x)}.
$$
 Define c.p.c.~ maps
$\zeta_{\rho},\zeta_{\sigma} \colon \Dh \to A$ by
$$
 \zeta_{\rho}(d) := \varphi_{\rho}(u_{\rho}(\be_{\Ch([0,1])} \otimes
 d\otimes \be_{\Dh}) u_{\rho}^*) \; , \;
 \zeta_{\sigma}(d) := \varphi_{\sigma}(u_{\sigma}(\be_{\Ch([0,1])} \otimes
d\otimes \be_{\Dh}) u_{\sigma}^*) \,.
$$
 Finally, we define
$ \psi \colon \Dh  \to A $ by
$$
\psi(d) := h_1 \cdot \rho(d)
        + h_2\cdot \zeta_{\rho}(d)
        + h_3\cdot \zeta_{\sigma}(d)
        + h_4 \cdot \sigma(d).
$$
$\psi$ is clearly a c.p.c.~ map. Note that
$\psi(d)_{[0,s-3\delta]} = \rho(d)_{[0,s-3\delta]}$ and
$\psi(d)_{[s+3\delta,1]} = \sigma(d)_{[s+3\delta,1]}$. In
particular, it follows that $\psi$ is $(\Fh;
\Gh',\varepsilon')$-good for some neighborhood of the endpoints of
the interval $[r,t]$, and is $(\Fh;\Gh,\varepsilon)$-good for
$[r,s-3\delta] \cup [s+3\delta,t]$.

It remains to show that $\psi$ is $(\Fh; \Gh,\varepsilon)$-good
for $[s-3\delta,s+3\delta]$. To verify condition (U) in Definition
\ref{good-maps-def}, note that
$$
 \psi(\be_{\Dh}) = h_1 \cdot
\rho(\be_{\Dh}) + h_2 \cdot \mu_{\rho}(\be_{\Dh} \otimes
\be_{\Dh}) + h_3 \cdot \mu_{\sigma}(\be_{\Dh}\otimes \be_{\Dh}) +
h_4\cdot \sigma(\be_{\Dh}).
 $$
Thus, for any $b \in \Fh$ and any $x \in [s-3\delta,s+3\delta]$ we
have
\begin{eqnarray*}
\lefteqn{ \|(b \psi(\be_{\Dh})  - b)_x\| }
    \\  & \mbox{}\leq &
    h_1(x)\underbrace{\|(b \rho(\be_{\Dh}) - b)_x\|}_{\mbox{}<\varepsilon'}
    \\ & & \mbox{} + h_2(x) \cdot  (
   \underbrace{\|(b  (\mu_{\rho}(\be_{\Dh} \otimes \be_{\Dh}) -
   \rho(\be_{\Dh})\nu_{\rho}(\be_{\Dh})))_x\|}_{\mbox{}
       <\varepsilon' \; \textrm{ by \ref{mu-lemma}
       (ii)}_{\rho}}
    +
     \|(b  \rho(\be_{\Dh})\nu_{\rho}(\be_{\Dh}) - b)_x\|  )
    \\ & & \mbox{}+
     h_3(x) \cdot (
     \underbrace{\|(b  (\mu_{\sigma}(\be_{\Dh} \otimes \be_{\Dh}) -
   \sigma(\be_{\Dh})\nu_{\sigma}(\be_{\Dh})))_x\|}_{\mbox{}
       <\varepsilon' \; \textrm{ by \ref{mu-lemma}
       (ii)}_{\sigma}} +
     \|(b  \sigma(\be_{\Dh})\nu_{\sigma}(\be_{\Dh}) - b)_x\|  )
    \\ & & \mbox{} +
    h_4(x)\underbrace{\|(b \sigma(\be_{\Dh}) -
    b)_x\|}_{\mbox{}<\varepsilon'}
  \\  & \mbox{}<  &
  \underbrace{(h_1(x) + h_2(x)+h_3(x)+h_4(x))}_{=1}\varepsilon'
    \\ & & \mbox{}+ h_2(x)\cdot \|(b  \rho(\be_{\Dh})\nu_{\rho}(\be_{\Dh}) - b)_x\|
    +
     h_3(x) \cdot \|(b \sigma(\be_{\Dh})\nu_{\sigma}(\be_{\Dh}) - b)_x\|
  \\ & \mbox{}\leq &
     h_2(x)\cdot  (
    \underbrace{\|(b  \rho(\be_{\Dh}) - b)_x\|}_{\mbox{}<\varepsilon'}
    \|\nu_{\rho}(\be_{\Dh})_x\|
     +
     \underbrace{\|(b\nu_{\rho}(\be_{\Dh}) - b)_x\|}_{\mbox{}<3\varepsilon'})
    \\ & & \mbox{}+
     h_3(x)\cdot  (
    \underbrace{\|(b  \sigma(\be_{\Dh}) - b)_x\|}_{\mbox{}<\varepsilon'}
    \|\nu_{\sigma}(\be_{\Dh})_x\|
     +
     \underbrace{\|(b\nu_{\sigma}(\be_{\Dh}) - b)_x\|}_{\mbox{}<3\varepsilon'})
    + \varepsilon'
  \\ &\mbox{}<& (h_2(x)+h_3(x))\cdot 4\varepsilon' + \varepsilon'
  \leq 5\varepsilon' < \varepsilon.
\end{eqnarray*}
For condition (C), let $b \in \Fh$, $d \in \Gh$, $x \in
[s-3\delta,s+3\delta]$. We wish to show that $\| b_x\psi(d)_x
-\psi(d)_xb_x \| < \varepsilon$. By the definition of $\psi$, and
the fact that $\rho,\sigma$ are $(\Fh;\Gh,\varepsilon)$-good for
$\{x\}$, it will suffice to show that
$$
\| b_x\zeta_{\rho}(d)_x -\zeta_{\rho}(d)_xb_x \| < \varepsilon
\mbox{ and }
 \| b_x\zeta_{\sigma}(d)_x -\zeta_{\sigma}(d)_xb_x \| < \varepsilon
$$
We will show it for $\zeta_{\rho}$ -- the proof for
$\zeta_{\sigma}$ is similar. If $x \in [s-3\delta,s-2\delta]$ then
$\zeta_{\rho}(d)_x = \mu_{\rho}(d)_x$, and hence
$\|b_x(\zeta_{\rho}(d)_x - \rho(d)_x\nu_{\rho}(d)_x)\| <
\varepsilon'$. Thus, we have
\begin{eqnarray*}
\lefteqn{\| b_x\zeta_{\rho}(d)_x -\zeta_{\rho}(d)_xb_x \| <
 \|[ b_x , \rho(d)_x\nu_{\rho}(d)_x]\| +
 2\varepsilon'}
 \\ & \leq &
 \underbrace{\|[ b_x , \rho(d)_x]\|}_{\mbox{}< \varepsilon'}
 +
 \underbrace{\|[ b_x , \nu_{\rho}(d)_x]\|}_{\mbox{}< 3\varepsilon'}
 + 2\varepsilon'
 < 6\varepsilon' < \varepsilon.
\end{eqnarray*}
If $x \in [s-\delta,s+3\delta]$, we have (by (\ref{u_1}))
$$\|\zeta_{\rho}(d)_x -
\mu_{\rho}(\be_{\Dh} \otimes d)_x\|< \frac{\varepsilon}{4}$$
 and by \ref{mu-lemma} (ii)$_{\rho}$,
$$\|b_x(\mu_{\rho}(\be_{\Dh} \otimes d)_x -
\rho(\be_{\Dh})_x\nu_{\rho}(d)_x)\| \, ,\, \|(\mu_{\rho}(\be_{\Dh}
\otimes d)_x - \rho(\be_{\Dh})_x\nu_{\rho}(d)_x)b_x\| <
\varepsilon'.
 $$
So, as in the previous consideration, we have
\begin{eqnarray*}
\| b_x\zeta_{\rho}(d)_x -\zeta_{\rho}(d)_xb_x \|
 <
\|[b_x,\rho(\be_{\Dh})_x\nu_{\rho}(d)_x]\| +
\frac{2\varepsilon}{4} + 2\varepsilon' < \frac{2\varepsilon}{4} +
6\varepsilon' \leq \varepsilon.
\end{eqnarray*}
Finally, if $x \in [s-2\delta,s-\delta]$, find $i$ such that
\begin{equation} \label{eq-z0}
\|u_{\rho\, x} - y_i \|< \frac{\varepsilon}{9} \,.
\end{equation}
 We then have
\begin{equation} \label{eq-z1}
 \|\zeta_{\rho}(d)_x - \mu_{\rho}(y_i(d \otimes
\be_{\Dh})y_i^*)_x\| <  \frac{2\varepsilon}{9} \,.
\end{equation}
Note that
\begin{equation} \label{eq-z2}
 \mu_{\rho}(y_i(d \otimes \be_{\Dh})y_i^*)
 = \sum_{j,k=1}^m \mu_{\rho}(v_{ji}dv_{ki}^* \otimes
 w_{ji}w_{ki}^*)
\end{equation}
and for any $1\leq j,k\leq m$, we have (by \ref{mu-lemma}
(ii)$_{\rho}$)
\begin{equation} \label{eq-z3}
 \|b_x(\mu_{\rho}(v_{ji}dv_{ki}^* \otimes
 w_{ji}w_{ki}^*) - \rho(v_{ji}dv_{ki}^*)
 \nu_{\rho}(w_{ji}w_{ki}^*))_x\| < \varepsilon'
\end{equation}
 and we get a similar estimate by placing $b_x$ on
the right rather than the left. Thus,
\begin{eqnarray*}
\lefteqn{\| b_x\zeta_{\rho}(d)_x -\zeta_{\rho}(d)_xb_x \|}\\
& \mbox{}< & \|[b_x, \mu_{\rho}(y_i(d \otimes
\be_{\Dh})y_i^*)_x]\| +  \frac{4\varepsilon}{9}
 \\ &\mbox{}< &
 \sum_{j,k=1}^m \|[b_x,\rho(v_{ji}dv_{ki}^*)_x
 \nu_{\rho}(w_{ji}w_{ki}^*)_x]\| +  \frac{4\varepsilon}{9} + 2m^2\varepsilon'
 \\ &\mbox{}\leq &
 \sum_{j,k=1}^m ( \underbrace{\|[b_x,\rho(v_{ji}dv_{ki}^*)_x]\|}_{\mbox{}<\varepsilon'} +
 \underbrace{\|[b_x,\nu_{\rho}(w_{ji}w_{ki}^*)_x]\|}_{\mbox{}<3\varepsilon'}  ) +  \frac{4\varepsilon}{9} + 2m^2\varepsilon'
 \\ &\mbox{}<&
 6m^2\varepsilon' +  \frac{4\varepsilon}{9} \leq \varepsilon.
\end{eqnarray*}

For condition (M), let $b \in \Fh$, $d,d' \in \Gh$, $x \in
[s-3\delta,s+3\delta]$. We must show that
$$\left \|(b(\psi(dd') - \psi(d) \psi(d')))_{x} \right \|<
\varepsilon.$$

We shall just show it for $x \in [s-3\delta,s]$ -- the case $x \in
[s,s+3\delta]$ is obtained similarly, where the roles of $\rho$
and $\sigma$ are reversed.

If $x \in [s-3\delta,s-2\delta]$,  We have that
$$
\psi(c)_x = h_1(x) \rho(c)_x + h_2(x) \mu_{\rho}(c \otimes
\be_{\Dh})_x.
$$
For any $c \in \Gh'$, $b \in \Fh$, we have (\ref{mu-lemma}
(ii)$_{\rho}$)
$$
 \|b_x(\mu_{\rho}(c
\otimes \be_{\Dh})_x -
 \rho(c)_x\nu_{\rho}(\be_{\Dh})_x)\| < \varepsilon'
$$
and furthermore
\begin{eqnarray*}
\lefteqn{\|b_x(\rho(c)_x\nu_{\rho}(\be_{\Dh})_x - \rho(c)_x)\| }
 \\ & \mbox{}\leq &
 \underbrace{\|b_x\nu_{\rho}(\be_{\Dh})_x -
 b_x\|}_{<3\varepsilon'}
 \|\rho(c)_x\|
 +
 \underbrace{\|b_x[\rho(c)_x,\nu_{\rho}(\be_{\Dh})_x]\|}_{
   \mbox{}<2\varepsilon' \; \textrm{ by \ref{mu-lemma}
   (i)}_{\rho}}
 < 5\varepsilon'.
\end{eqnarray*}
Thus, we see that
\begin{equation} \label{mu-psi-rho}
\|b_x(\mu_{\rho}(c\otimes \be_{\Dh})_x - \rho(c)_x)\| \; , \;
 \|b_x(\psi(c)_x -
\rho(c)_x)\|<6\varepsilon'.
\end{equation}
So,
\begin{eqnarray*}
\lefteqn{\|(b(\psi(dd') - \psi(d) \psi(d')))_{x}  \|}
 \\ &\mbox{}&<
 \|(b(\rho(dd') - \rho(d)\mu_{\rho}(d' \otimes \be_{\Dh})))_x\| +
 12\varepsilon'
  \\ & \mbox{}&\leq
 \underbrace{\|(b(\rho(dd') - \rho(d)\rho(d')))_x\|}_{<\varepsilon'} +
 \|b_x\rho(d)_x(\rho(d') - \mu_{\rho}(d' \otimes \be_{\Dh}))\| + 12\varepsilon'
 \\ &\mbox{}&<
  \underbrace{\|\rho(d)_xb_x(\rho(d') - \mu_{\rho}(d' \otimes
  \be_{\Dh}))\|}_{\mbox{}<6\varepsilon' \; \textrm{ by
  (\ref{mu-psi-rho})}}
  + \underbrace{\|[b_x,\rho(d)_x]\|}_{\mbox{}<\varepsilon'}
  + 13\varepsilon' \\
  & \mbox{}&< 20\varepsilon' < \varepsilon.
\end{eqnarray*}

If $x \in [s-2\delta,s-\delta]$, we have
 $\psi(c)_x =  \zeta_{\rho}(c)_x = \mu_{\rho}(u_{\rho \, x}(c \otimes \be_{\Dh})u_{\rho \, x}^*)$.
Repeating the considerations from equations
(\ref{eq-z0})--(\ref{eq-z3}) above, we find an $i$ such that
\begin{equation} \label{M-uy}
\|u_{\rho\, x} - y_i \|< \frac{\varepsilon}{9}
\end{equation}
and
then we have
\begin{eqnarray*}
\lefteqn{\left \|(b(\psi(dd') - \psi(d) \psi(d')))_{x} \right \|
 }
 \\ &\mbox{} < &
 \left \| b_x \left (  \zeta_{\rho}(dd') -
  \!\!\!\! \sum_{j,k,\ell,n=1}^m \!\!\!\! \mu_{\rho}(v_{ji}dv_{ki}^* \otimes
 w_{ji}w_{ki}^*) \mu_{\rho}(v_{\ell i}d'v_{ni}^* \otimes
 w_{\ell i}w_{ni}^*) \right)_x \right \| + \frac{4\varepsilon}{9}
 \,.
\end{eqnarray*}
Now, note that if $c_1,c_2,c_3,c_4 \in \Gh'$ (which we will apply
as
 $c_1 = v_{ji}dv_{ki}^*$, $c_2 = w_{ji}w_{ki}^*$,
 $c_3 = v_{\ell i}d'v_{ni}^*$, $c_4 =w_{\ell i}w_{ni}^*$), we
 have (by \ref{mu-lemma} (ii)$_{\rho}$)
$$\|b_x(\mu_{\rho}(c_1 \otimes
 c_2)_x -
 \rho(c_1)_x\nu_{\rho}(c_2)_x)\| <
 \varepsilon' $$
and therefore
\begin{eqnarray}
\label{rho-nu1}
 \lefteqn{\|b_x(\mu_{\rho}(c_1 \otimes
 c_2)_x\mu_{\rho}(c_3 \otimes
 c_4)_x -
 \rho(c_1)_x\nu_{\rho}(c_2)_x
 \rho(c_3)_x\nu_{\rho}(c_4)_x)\| }
 \\\nonumber & < &
  \|b_x \rho(c_1)_x\nu_{\rho}(c_2)_x
  ( \mu_{\rho}(c_3 \otimes
 c_4)_x - \rho(c_3)_x\nu_{\rho}(c_4)_x ) \| + \varepsilon'
 \\\nonumber & < &
 \|b_x( \mu_{\rho}(c_3 \otimes
 c_4)_x - \rho(c_3)_x\nu_{\rho}(c_4)_x )\| +
  \|[b_x,\rho(c_1)_x\nu_{\rho}(c_2)_x]\| +
  \varepsilon'
 \\\nonumber & < &
  \underbrace{\|[b_x,\rho(c_1)_x]\|}_{<\varepsilon'} +
  \underbrace{\|[b_x,\nu_{\rho}(c_2)_x]\|}_{<3\varepsilon'}
   + 2\varepsilon' <
  6\varepsilon'.
\end{eqnarray}
Thus, we have that
\begin{eqnarray} \label{rho-nu2}
\lefteqn{\left \|(b(\psi(dd') - \psi(d) \psi(d')))_{x} \right \|
 }
  \\\nonumber  &< &
 \left \| b_x \! \left ( \! \zeta_{\rho}(dd') -
  \!\!\!\!   \sum_{j,k,\ell,n=1}^m \!\!\!\! \rho(v_{ji}dv_{ki}^*)
 \nu_{\rho}(w_{ji}w_{ki}^*) \rho(v_{\ell i}d'v_{ni}^*)
 \nu_{\rho}(w_{\ell i}w_{ni}^*) \!\right)_x \right \|
 \\\nonumber & & \mbox{}+\frac{4\varepsilon}{9} +
 6m^4\varepsilon'.
\end{eqnarray}
Now, note that
\begin{eqnarray*} \label{rho-nu3}
\lefteqn{\|b_x(\rho(c_1)
 \nu_{\rho}(c_2) \rho(c_3)
 \nu_{\rho}(c_4) -
 \rho(c_1) \rho(c_3)
 \nu_{\rho}(c_2)
 \nu_{\rho}(c_4))_x\| }
 \\\nonumber & \leq &
 \underbrace{\|[b_x,\rho(c_1)_x]\|}_{<\varepsilon'} +
 \underbrace{\|b_x[\nu_{\rho}(c_2)_x ,
 \rho(c_3)_x]\|}_{<2\varepsilon'}
 < 3\varepsilon'
\end{eqnarray*}
and
\begin{eqnarray*} \label{rho-nu4}
\lefteqn{  \|b_x(
 \rho(c_1) \rho(c_3)
 \nu_{\rho}(c_2)
 \nu_{\rho}(c_4)
 -
 \rho(c_1c_3)
 \nu_{\rho}(c_2c_4)
 )_x\|
 }
 \\\nonumber & \leq &
 \|b_x\rho(c_1)_x \rho(c_3)_x (\nu_{\rho}(c_2)
 \nu_{\rho}(c_4) - \nu_{\rho}(c_2c_4) )_x\|
 +
 \underbrace{\|b_x(\rho(c_1) \rho(c_3) -
 \rho(c_1c_3))_x\|}_{<\varepsilon'}
 \\\nonumber & < &
 \|[b_x,\rho(c_1)_x \rho(c_3)_x]\| +
 \underbrace{\|b_x(\nu_{\rho}(c_4) - \nu_{\rho}(c_2c_4)
 )_x\|}_{<3\varepsilon'}
 + \varepsilon'
 \\\nonumber & < &
 \|[b_x,\rho(c_1)_x]\| + \|[b_x,\rho(c_3)_x]\| + 4\varepsilon'
 < 6\varepsilon'
\end{eqnarray*}
and so
\begin{equation} \label{rho-nu5}
\|b_x(\rho(c_1)
 \nu_{\rho}(c_2) \rho(c_3)
 \nu_{\rho}(c_4) -
 \rho(c_1c_3)
 \nu_{\rho}(c_2c_4)
 )_x\| < 9\varepsilon'.
\end{equation}
 Thus,
\begin{eqnarray*}
\lefteqn{\left \|(b(\psi(dd') - \psi(d) \psi(d')))_{x} \right \|
 }
  \\ & \stackrel{(\ref{rho-nu2},\ref{rho-nu5})}{<} &
 \left \| b_x \left (  \zeta_{\rho}(dd') -
  \!\!\!\!   \sum_{j,k,\ell,n=1}^m \!\!\!\!
  \rho(v_{ji}dv_{ki}^*v_{\ell i}d'v_{ni}^*)
 \nu_{\rho}(w_{ji}w_{ki}^*w_{\ell i}w_{ni}^*)
 \right)_x \right \|\\
 & &  + \frac{4\varepsilon}{9} +
 15m^4\varepsilon'
 \\ & \stackrel{\ref{mu-lemma} \textrm{ (ii)}_{\rho}}{<} &
 \left \| b_x \left (  \zeta_{\rho}(dd') -
  \!\!\!\!   \sum_{j,k,\ell,n=1}^m \!\!\!\!
  \mu_{\rho}(v_{ji}dv_{ki}^*v_{\ell i}d'v_{ni}^*
  \otimes w_{ji}w_{ki}^*w_{\ell i}w_{ni}^*)
 \right)_x \right \| \\
 & & + \frac{4\varepsilon}{9} +
 16m^4\varepsilon'
 \\ & = &
 \left \| b_x \left (  \zeta_{\rho}(dd') -
 \mu_{\rho}(y_i(d \otimes \be_{\Dh}) y_i^*
  y_i(d' \otimes \be_{\Dh}) y_i^*) \right )_x \right \|
 + \frac{4\varepsilon}{9} +
 16m^4\varepsilon'
 \\ & \stackrel{(\ref{M-uy})}{<} &
 \ \| b_x  (  \zeta_{\rho}(dd') -
 \underbrace{\mu_{\rho}(u_{\rho \, x}(d \otimes \be_{\Dh})
  u_{\rho \, x}^* u_{\rho \, x}
  (d' \otimes \be_{\Dh}) u_{\rho \, x}^*)}_{=\zeta_{\rho}(dd')}
    )_x \| \\
    & &
 +\frac{8\varepsilon}{9} +
 16m^4\varepsilon'
 \\ & = &  \frac{8\varepsilon}{9} + 16m^4\varepsilon' = \varepsilon.
\end{eqnarray*}
Finally, if $x \in [s-\delta,s]$, we have that
\begin{align*} \lefteqn{\psi(c)_x =
h_2(x)\zeta_{\rho}(c)_x + h_3(x)\zeta_{\sigma}(c)_x}
 \\ & =
 h_2(x)\mu_{\rho}(u_1 (c\otimes \be_{\Dh}) u_1^*)_x + h_3(x)
 \mu_{\sigma}(u_1 (c\otimes \be_{\Dh}) u_1^*)_x \,.
\end{align*}
Note that if $c \in \Gh \cdot \Gh \subseteq \Gh'$, we have (by
(\ref{u_1}))
\begin{eqnarray*}
\|\mu_{\rho}(u_1 (c\otimes \be_{\Dh}) u_1^*) -
 \mu_{\rho}(\be_{\Dh} \otimes c)\| < \frac{\varepsilon}{4}
\end{eqnarray*}
and if $b \in \Fh$, then we have (using \ref{mu-lemma}
(ii)$_{\rho}$ and condition (U) for $\rho$)
\begin{eqnarray*}
\lefteqn{\|b_x(\mu_{\rho}(\be_{\Dh} \otimes c)_x -
\nu_{\rho}(c)_x)\| }
 \\ & \leq &
 \|b_x(\mu_{\rho}(\be_{\Dh} \otimes c)_x -\rho (\be_{\Dh})_x
\nu_{\rho}(c)_x)\|
 +
 \|(b_x\rho (\be_{\Dh})_x - b_x)\nu_{\rho}(c)_x\| \\
 & < & 2\varepsilon'.
\end{eqnarray*}
Thus,
$$
\|b_x(\zeta_{\rho}(c)_x - \nu_{\rho}(c)_x)\| <
\frac{\varepsilon}{4} + 2\varepsilon' \,.
$$
Similarly,
$$
\|b_x(\zeta_{\sigma}(c)_x - \nu_{\sigma}(c)_x)\| <
\frac{\varepsilon}{4} + 2\varepsilon' \,.
$$
But $\|b_x(\nu_{\rho}(c)_x - \nu_{\sigma}(c)_x\| < 2\varepsilon'$
(by \ref{mu-lemma} (iii)), and thus we have
$$
\|b_x(\zeta_{\sigma}(c)_x - \nu_{\rho}(c)_x)\| <
\frac{\varepsilon}{4} + 4\varepsilon' \,.
$$
Therefore,
$$
\|b_x(\psi(c)_x - \nu_{\rho}(c)_x)\| < \frac{\varepsilon}{4} +
4\varepsilon'.
$$
So, for $d,d' \in \Gh$, we have (since $\nu_{\rho}$ is
$(\Fh;\Gh',3\varepsilon')$-good for $\{x\}$)
\begin{eqnarray*}
\lefteqn{
 \|b_x ( \psi(dd') - \psi(d)\psi(d') )_x\|  }
 \\ & < &
  \|b_x ( \nu_{\rho}(dd') - \nu_{\rho}(d)\psi(d') )_x\|
  + \frac{2\varepsilon}{4} + 8 \varepsilon'
 \\ & \leq &
  \underbrace{\|b_x ( \nu_{\rho}(dd') - \nu_{\rho}(d)\nu_{\rho}(d') )_x\|}_{<3\varepsilon'}
  + \underbrace{\|[b_x,\nu_{\rho}(d)_x]\|}_{<3\varepsilon'}
  \\ & & \mbox{}
  + \|\nu_{\rho}(d)\| \underbrace{\|b_x(\psi(d')_x -
  \nu_{\rho}(d')_x)\|}_{<\,\varepsilon/4 +
4\varepsilon'}
  +  \frac{2\varepsilon}{4} + 8 \varepsilon' \\
  &
  < &
  \frac{3\varepsilon}{4} + 18\varepsilon' <\varepsilon.
\end{eqnarray*}
\end{nproof}
\en

\bn
\label{D-stable-bundles}
\begin{thms}
Let $\Dh$ be a $K_{1}$-injective strongly self-absorbing
$C^{*}$-algebra. Let $X$ be a locally compact metrizable space of
finite covering dimension. Let $A$ be a separable
$\Ch_{0}(X)$-algebra. It follows that $A$ is $\Dh$-stable if and
only if $A_{x}$ is $\Dh$-stable for each $x \in X$.
\end{thms}

\begin{nproof}
If $A$ is $\Dh$-stable, then so are all fibres $A_{x}$, since
$\Dh$-stability passes to quotients by
\cite[Corollary~3.3]{TomsWinter:ssa}. We prove the converse. By
\ref{C(X)-limits}, it suffices to restrict to compact $X$. By
\cite[Theorem~V.3]{HurWall:Dim} we may assume that $X$ is a subset
of $[0,1]^{N}$ for some $N \in \N$. By \ref{restriction comment},
we can furthermore assume that $X = [0,1]^{N}$ (the fibres are
either the original fibres, or 0, all of which are $\Dh$-stable,
so the hypothesis holds).

We can now proceed by induction. For $N=1$, we use Proposition
\ref{D-stability-characterization} b). Given finite sets $\Fh
\subset A$ and $\be_{\Dh} \in \Gh \subset \Dh$, and
$\varepsilon>0$, we would like to find  a c.p.c.~ map $\psi\colon
\Dh \to A$ which is $(\Fh ; \Gh , \varepsilon)$-good for $[0,1]$.
We may assume without loss of generality that $\Fh,\Gh$ are
self-adjoint and consist of elements of norm at most one.

Select $\Gh' \supset \Gh$, $\varepsilon>\varepsilon'>0$ as in
Lemma \ref{patching lemma}. Use Lemma \ref{commuting-good-maps} to
find some natural number $n$, points $0 = t_0 < t_1 < \cdots < t_n
= 1$ and c.p.c.~ maps $\psi_k \colon \Dh \to A$ such that $\psi_k$
is $(\Fh ; \Gh',\varepsilon')$-good for $[t_{k-1},t_k]$. Note that
any map which is $(\Fh;\Gh',\varepsilon')$-good for some interval
is trivially also $(\Fh ; \Gh,\varepsilon ;
\Gh',\varepsilon')$-good for the same interval. Now, by a repeated
application ($n-1$ times) of Lemma \ref{patching lemma}, we can
find a c.p.c.~ map $\psi\colon \Dh \to A$ which is $(\Fh ;
\Gh,\varepsilon ; \Gh',\varepsilon')$-good for $[0,1]$. Such a map
is in particular $(\Fh ; \Gh,\varepsilon)$-good for $[0,1]$, i.e.~
satisfies the required conditions.

For the induction step, suppose now that $N>1$, and that the
statement is true for $X = [0,1]^k$ for all $k<N$. We may clearly
regard $A$ as a $\Ch([0,1])$-algebra with fibres $A_{\{t\} \times
[0,1]^{N-1}}$ for $t \in [0,1]$. By the base case of the
induction, if $A_{\{t\} \times [0,1]^{N-1}}$ is $\Dh$-stable for
all $t \in [0,1]$, then so is A. But, for each $t \in [0,1]$, and
each $x \in \{t\} \times [0,1]^{N-1}$, $(A_{\{t\} \times
[0,1]^{N-1}})_{x}=A_{x}$ is $\Dh$-stable by assumption. Upon
identifying $\{t\}\times [0,1]^{N-1}$ with $[0,1]^{N-1}$, it
follows from the induction hypothesis that $A_{\{t\} \times
[0,1]^{N-1}}$ also absorbs $\Dh$.
\end{nproof}
\en

The two examples below show that one cannot drop the finite dimensionality
condition from the preceding result.

\bn
\begin{exs}
\label{UHFexample} We give here an example of a $\Ch(X)$-algebra
$A$, where $X = \prod_{n=1}^\infty S^2$ is infinite dimensional,
whose fibres $A_x$ are isomorphic to the CAR-algebra (the
UHF-algebra of type $2^\infty$) for all $x \in X$, and such that
$A$ does not absorb this UHF-algebra. Actually, one cannot embed
$M_2$ unitally into $A$.

The $K_0$-group of $\Ch(S^2)$ is isomorphic to $\Z \oplus \Z$ and
is generated (as a $\Z$-module) by the $K_0$-classes of two
1-dimensional projections $e$ and $p$, where $e$ is trivial (=
constant) and where $p$ is the ``Bott projection''. The latter
projection belongs to $M_2(\Ch(S^2))$. Hence we can choose $e$ and
$p$ to be mutually orthogonal projections in $M_3(\Ch(S^2))$. Put
$$B = (e+p)M_3(\Ch(S^2))(e+p), \qquad A = \bigotimes_{n=1}^\infty
B.$$

Note that $B$ is a $\Ch(S^2)$-algebra with fibres $B_x = M_2$ (the
projections $e$ and $p$ are not equivalent in $B$, but $e_x$ is
equivalent to $p_x$ for all $x \in S^2$). It follows from
Lemma~\ref{lm:inftensor} that $A$ is a $\Ch(X)$-algebra with fibres $A_x$
isomorphic to the CAR-algebra $\bigotimes_{n=1}^\infty M_2$ for all $x
\in X$.

We proceed to show that one cannot embed $M_2$ unitally into $A$. We
do so by showing that $[1_A] \in K_0(A)$ is not divisible by $2$. By
continuity of $K_0$ it suffices to show that the class in $K_0$ of the
unit, call it $q_m$, of $\bigotimes_{n=1}^m B$ is not divisible by
$2$. Now,
\begin{eqnarray*}
q_m & = & (e+p) \otimes (e+p) \otimes (e+p) \otimes \cdots \otimes
(e+p) \\&=&  \sum_{I \subseteq \{1,2,\dots,m\}} p_I,
\end{eqnarray*}
where
$$p_I = r_1 \otimes r_2 \otimes \cdots \otimes r_m, \qquad r_j =
\begin{cases} p, & j \in I \\ e, & j \notin I \end{cases}.$$

We already noted that $[e]$ and $[p]$ form a basis for $K_0(\Ch(S^2))$
as a $\Z$-module. It follows from the K\"unneth theorem (and by
induction) that the $2^m$ elements $[p_I]$, $I \subseteq
\{1,2,\dots,m\}$, form a basis for $K_0(\Ch((S^2)^m))$ (again as a
$\Z$-module). This shows that $[q_m] = \sum_I [p_I]$ is not divisible
by $2$ in $K_0(\Ch((S^2)^m))$.

This example can easily be amended to give, for each UHF-algebra $\Dh$
of infinite type, a $\Ch(X)$-algebra $A$ (with $X$ as above) with fibres
$A_x$ isomorphic to $\Dh$ for all $x \in X$ and such that $A$ does not
absorb $\Dh$ tensorially. Indeed, one can construct $A$ such that one
cannot embed any non-trivial matrix algebra unitally into $A$.
\end{exs}
\en

\bn
\begin{exs}
\label{Zexample} We give here an example of a $\Ch(X)$-algebra
$A$, where, as in the previous example, $X = \prod_{n=1}^\infty
S^2$ and the fibres $A_x$ are all isomorphic to the CAR-algebra
(the UHF-algebra of type $2^\infty$)---in particular, the fibres
absorb the Jiang--Su algebra $\Zh$ (see \cite{JiaSu:Z})---but such
that $A$ does not absorb $\Zh$. In fact, we show that the
semigroup $V(A)$ of Murray-von Neumann equivalence classes of
projections in matrix algebras over $A$ is not almost
unperforated, i.e., there are elements $x,y \in V(A)$ and a
natural number $n$ such that $(n+1)x \le ny$ but $x \nleq y$. It
then follows from \cite[Corollary~4.8]{Ror:Z-absorbing} that $A$
cannot be $\Zh$-absorbing.

Here, as in Example~\ref{UHFexample} above, we also get a
$\Ch(X)$-algebra with fibres isomorphic to the CAR-algebra, but
which itself does not absorb the CAR-algebra. (Any \Cs{} that
tensorially absorbs a UHF-algebra will also absorb the Jiang--Su
algebra, because any UHF-algebra absorbs the Jiang--Su algebra,
cf.\ \cite{JiaSu:Z}.) All the same, we included
Example~\ref{UHFexample}, as it is technically easier than the
present example.

As in Example~\ref{UHFexample}, let $p \in M_2(\Ch(S^2))$ be the
one-dimensional ``Bott projection''. For each natural number $m$
identify the two \Cs s
$\bigotimes_{n=1}^m M_2(\Ch(S^2))$ and $M_{2^m}(\Ch((S^2)^m))$,
and find in $M_{2^m+1}(\Ch((S^2)^m))$ mutually orthogonal projections
$e$ and $p^{\otimes m}$, such that $e$ is a trivial one-dimensional projection
and $p^{\otimes m}$ is (equivalent to)
$$p \otimes p \otimes \cdots \otimes p \in
M_{2^m}(\Ch((S^2)^m)).$$

Put
$$m(1)=m(2)=1, \qquad m(j) = 2^{j-2}, \quad j \ge 3,$$
and put
$$B_j = (e+p^{\otimes
  m(j)})M_{2^{m(j)}+1}(\Ch((S^2)^{m(j)}))(e+p^{\otimes m(j)}),
\qquad A = \bigotimes_{j=1}^\infty B_j.$$
As in Example~\ref{UHFexample} (again using Lemma~\ref{lm:inftensor})
we see that $A$ is a $\Ch(X)$-algebra with fibres $A_x$ isomorphic to
the CAR-algebra for all $x \in X$.

For each $n \ge 2$ consider the two projections
\begin{eqnarray*}
f_n &=& e \otimes e \otimes 1_{B_3} \otimes \cdots \otimes 1_{B_n}\\
&=& e \otimes e \otimes (e+p^{\otimes m(3)}) \otimes \cdots \otimes
(e+p^{\otimes m(n)}),\\
g_n &=& (p \otimes e + e \otimes p) \otimes 1_{B_3} \otimes \cdots
\otimes 1_{B_n}\\
&=& (p \otimes e + e \otimes p) \otimes (e+p^{\otimes m(3)}) \otimes
\cdots \otimes (e+p^{\otimes m(n)})
\end{eqnarray*}
in $ \bigotimes_{j=1}^n B_j$. Let $f$ and $g$ be the corresponding
projections in $A$. We show below that $4 [f] \le 3[g]$ and that $[f]
\nleq [g]$ in $V(A)$. This will settle the
claims made in the first paragraph of this example.

The projection $f_2 \oplus f_2 \oplus f_2 \oplus f_2$ has dimension
$4$, the projection $g_2 \oplus g_2 \oplus g_2$ has dimension $6$,
and as the difference $6-4=2$ is greater than or equal to
$3/2=(\dim((S^2)^2)-1)/2$ it follows from
\cite[9.1.2]{Hus:fibre}
that $f_2 \oplus f_2 \oplus f_2 \oplus f_2$ is equivalent to a
subprojection of $g_2
\oplus g_2 \oplus g_2$. This entails that $f \oplus f \oplus f \oplus
f \precsim g \oplus g \oplus g$, whence $4[f] \le 3[g]$ in $V(A)$.

To show that $f \nprecsim g$ it suffices to show that $f_n \nprecsim
g_n$ for all $n \ge 2$. We show that the Euler class of $g_n$ is
non-zero for all $n$. Hence $g_n$ cannot dominate any trivial
projection, and as $f_n$ does dominate a trivial projection, $f_n$ is
not equivalent to a subprojection of $g_n$.

To calculate the Euler class of $g_n$ let us first note that $g_n$
belongs to a matrix algebra over
$$\bigotimes_{j=1}^n \Ch((S^2)^{m(j)}) \; \cong \; \bigotimes_{j=1}^{M(n)}
\Ch(S^2) \; \cong \; \Ch((S^2)^{M(n)}),$$
where $M(j) = m(1) + m(2) + \cdots + m(j)$. Expand $g_n$ as follows:
\begin{eqnarray*}
g_n &=& (p \otimes e + e \otimes p) \otimes (e+p^{\otimes m(3)}) \otimes
\cdots \otimes (e+p^{\otimes m(n)})\\
&=& \sum_{I \subseteq \{3,4,\dots, n\}} p \otimes e \otimes
  q_I \;  + \; \sum_{I \subseteq
    \{3,4,\dots, n\}} e \otimes p \otimes q_I
\end{eqnarray*}
where
$$q_I = r_3 \otimes r_4 \otimes \cdots \otimes r_n, \qquad r_i =
\begin{cases} p^{\otimes m(i)}, & i \in I \\ e, & i \notin I. \end{cases}$$
Following the notation of \cite{Ror:simple}, for each subset $J$
of $\{1,2,\dots, M(n)\}$ put
$$p_J = r_1 \otimes r_2 \otimes \cdots \otimes r_{M(n)}, \qquad r_j =
\begin{cases} p, & j \in J,\\ e, & j \notin J, \end{cases}$$
which belongs to a matrix algebra over $\Ch((S^2)^{M(n)}) \cong
\bigotimes_{j=1}^{M(n)}
\Ch(S^2)$.

For $3 \le j \le n$ let $H_j$ be the set of integers $i$
such that $M(j-1)+1 \le i \le M(j)$,
and for $I \subseteq  \{3,4,\dots, n\}$ put
$$J_1(I) = \{1\} \cup \bigcup_{i \in I} H_i, \qquad J_2(I) = \{2\}
\cup \bigcup_{i \in I} H_i.$$
Then
$$p \otimes e \otimes q_I \sim p_{J_1(I)}, \qquad e \otimes p
\otimes q_I \sim p_{J_2(I)},$$
whence
$$g_n \sim \bigoplus_{I \subseteq \{3,4,\dots, n\}} p_{J_1(I)} \oplus
\bigoplus_{I \subseteq \{3,4,\dots, n\}} p_{J_2(I)}.
$$
We can now use \cite[Proposition 3.2 and Lemma 4.1]{Ror:simple} to
conclude that the Euler class of $g_n$ is non-zero. We have to
show that the family
$$\{J_1(I) \mid I \subseteq \{3,4,\dots,n\}\} \cup \{J_2(I) \mid I
\subseteq \{3,4,\dots,n\}\}$$
admits a matching (cf.\ \cite[Proposition 3.2 (iii)]{Ror:simple}). Now,
$J_1(\emptyset) = \{1\}$ and $J_2(\emptyset) = \{2\}$, so it is clear
how to match these two sets. If $I \subseteq \{3,4,\dots, n\}$ is
non-empty, then
choose the matching elements for $J_1(I)$ and for $J_2(I)$ in the set
$H_{\max I}$. This is possible because there are $2^{k -3}$ subsets
$I$ of $\{3,4,\dots, n\}$ with $\max I = k$ and there are $2 \cdot
2^{k-3} = m(k)$ elements in $H_k$.
\end{exs}
\en

\bn \label{pullbacks} We have just seen  that Theorem
\ref{D-stable-bundles}  does not remain true for infinite
dimensional $X$. However, we shall show in Proposition
\ref{locally D-stable} that when $\Dh$-stability of the fibres is
replaced by `local' $\Dh$-stability, a statement similar to
\ref{D-stable-bundles} holds. We first show that $\Dh$-stability
passes to pullbacks.

\begin{props}
Let $\Dh$ be a $K_{1}$-injective strongly self-absorbing
$C^{*}$-algebra. Suppose we have a pullback diagram of separable
$C^*$-algebras
$$\xymatrix{
 C \ar[r]^{\pi_1} \ar[d]_{\pi_2} & A_1 \ar[d]^{\varphi_1}
  \\
  A_2 \ar[r]_{\varphi_2} & B   }
$$
where $C$ is the pullback of $(A_1,A_2)$ along
$(\varphi_1,\varphi_2)$, and at least one of the maps
$\varphi_1,\varphi_2$ is surjective. If $A_1$ and $A_2$ are
$\Dh$-stable, then so is $C$.
\end{props}

\begin{nproof}
Let us assume that $\varphi_2$ is surjective. Therefore $\pi_1$ is
also surjective, and so, $C$ is an extension of $A_1$ by
$\ker(\pi_1)$. We may identify
$$C = \{(a_1,a_2) \in A_1\oplus A_2 \mid \varphi_1(a_1) =
\varphi_2(a_2)\} \subseteq A_1 \oplus A_2.$$ Under this
identification, we have $\ker(\pi_1) = \{(0,a) \mid \varphi_2(a) =
0\} = 0 \oplus \ker(\varphi_2)$. Since $A_2$ is $\Dh$-stable, and
$\Dh$-stability passes to ideals, it follows that
$\ker(\varphi_2)$, and therefore $\ker(\pi_1)$, are $\Dh$-stable.
Since $A_1$ and $\ker(\pi_1)$ are $\Dh$-stable, and
$\Dh$-stability passes to extensions by \cite[Theorem~4.3]{TomsWinter:ssa}, it follows that $C$ is
$\Dh$-stable, as required.
\end{nproof}
\en

\bn
\begin{rems}
In Proposition \ref{pullbacks}, we assumed that at least one of
the maps into $B$ is surjective. This assumption cannot be
removed. To see that, consider the following pullback diagram:
$$
\xymatrix{ \C \ar[r]^{\lambda \mapsto \lambda 1} \ar[d]_{\lambda \mapsto \lambda 1} & \Dh \ar[d]^{d \mapsto d \otimes 1} \\
\Dh \ar[r]_<<<<<{d \mapsto 1 \otimes d} &  \Dh \otimes \Dh}
$$
Even though the two copies of $\Dh$ are of course $\Dh$-stable,
$\C$ is not.
\end{rems}
\en

\bn \label{locally D-stable}
\begin{props}
Let $\Dh$ be a $K_{1}$-injective strongly self-absorbing
$C^{*}$-algebra. Let $X$ be a locally compact metrizable space
(not necessarily of finite dimension) and $A$ a separable
$\Ch_{0}(X)$-algebra. Suppose each $x \in X$ has a compact
neighborhood $V_{x}\subset X$ such that $A_{V_{x}}$ is
$\Dh$-stable for each $x \in X$. Then, $A$ is $\Dh$-stable.
\end{props}

\begin{nproof}
By \ref{C(X)-limits}  it suffices to prove the assertion for
compact $X$. But then we may assume that $X$ is covered by
finitely many closed subsets $V_{i}$, $i=1, \ldots,k$, such that
each $A_{V_{i}}$ is $\Dh$-stable. Set
\[
K_{\ell}:=\bigcup_{i=1}^{\ell} V_{i}.
\]
We prove that $A_{K_{\ell}}$ is $\Dh$-stable for each $\ell=1,
\ldots, k$ by induction; this will suffice as $X=K_{k}$.
$A_{K_{1}}=A_{V_{1}}$ is $\Dh$-stable by assumption. Suppose now
that $A_{K_{\ell-1}}$ is $\Dh$-stable for some $\ell \in \{2,
\ldots, k\}$. We may write  $A_{K_{\ell}}$ as a pullback
\[
\xymatrix{A_{K_{\ell}} \ar[r] \ar[d] & A_{K_{\ell-1}} \ar[d]  \\
A_{V_{\ell}} \ar[r] & A_{K_{\ell-1}\cap V_{\ell}}\, ,}
\]
where the maps are the restriction maps, and in particular are
surjective. By Proposition \ref{pullbacks}, we see that
$A_{K_{\ell}}$ is $\Dh$-stable, as required.
\end{nproof}
\en

\bn
\begin{rems}
Recall (\cite{BlaKumRor:apprdiv}) that a separable $C^*$-algebra
$A$ is said to be \emph{approximately divisible}  if there is a
central sequence of unital embeddings of $M_2 \oplus M_3$ into
$\Mh(A)$. Approximate divisibility is seen as a certain regularity
property, which is weaker than $\Dh$-stability for $\Dh$ a UHF
algebra, and stronger than $\Zh$-stability. One might ask if
similar results to the ones we have obtained can be found for
approximate divisibility. However, the following example shows
that the analogues of both  Theorem \ref{D-stable-bundles} and
Proposition \ref{locally D-stable} fail in this case.

Denote by $M_{2^{\infty}}, M_{3^{\infty}}$ the UHF algebras of
types $2^{\infty},3^{\infty}$ respectively. Let $$A = \{f \in
\Ch([0,1],M_{2^{\infty}} \otimes M_{3^{\infty}}) \mid  f(0) \in
M_{2^{\infty}} \otimes 1 \, , \, f(1) \in 1 \otimes M_{3^{\infty}}
\}$$ The embedding of $\Ch([0,1])$ into $A$ as the scalar
functions gives $A$ the structure of a $\Ch([0,1])$-algebra. The
fibres are $M_{2^{\infty}}$ at 0, $M_{3^{\infty}}$ at 1, and
$M_{2^{\infty}} \otimes M_{3^{\infty}}$ elsewhere, so in
particular, each fibre is approximately divisible. Furthermore,
$A_{[0,2/3]}$ and $A_{[1/3,1]}$ are approximately divisible (they
absorb $M_{2^{\infty}}$ and $M_{3^{\infty}}$, respectively). Thus,
$A$ satisfies analogous conditions to those of both Theorem
\ref{D-stable-bundles} and Proposition \ref{locally D-stable},
when one replaces $\Dh$-stability by approximate divisibility.
However, it is easy to see that $A$ is unital with no nontrivial
projections, and hence cannot be approximately divisible.
\end{rems}
\en

\bibliographystyle{amsplain}

\begin{thebibliography}{10}

\bibitem{Bla:encyc}
B.~Blackadar, \emph{{Operator Algebras}}, vol. 122, Encyclopaedia
of
  Mathematical Sciences. Subseries: Operator Algebras and Non-commutative
  Geometry, no. III, Springer Verlag, Berlin, Heidelberg, 2006.

\bibitem{BlaHan:quasitrace}
B.~Blackadar and D.~Handelman, \emph{{Dimension functions and
traces on
  $C^*$-algebras}}, J. Funct. Anal. \textbf{45} (1982), 297--340.

\bibitem{BlaKumRor:apprdiv}
B.~Blackadar, A.~Kumjian, and M.~R{\o}rdam, \emph{{Approximately
central matrix
  units and the structure of non-commutative tori}}, $K$-theory \textbf{6}
  (1992), 267--284.

\bibitem{Blanchard:C(X)-remarks}
E.~Blanchard, \emph{{A few remarks on exact $\Ch(X)$-algebras}},
Rev. Roum.
  Math. Pures Appl. \textbf{45} (2000), 565--576.

\bibitem{BlanchardKirchberg:Glimm}
E.~Blanchard and E.~Kirchberg, \emph{{ Global Glimm halving for
  $C^*$-bundles}}, J. Operator Theory \textbf{52} (2004), no.~2, 385--420.

\bibitem{BlanchardKirchberg:piHausdorff}
\bysame, \emph{{Non simple purely infinite $C^*$-algebras: The
Hausdorff
  case}}, J. Funct. Anal. \textbf{207} (2004), 461--513.

\bibitem{Cuntz:dimension}
J.~Cuntz, \emph{{Dimension functions on simple $C^*$-algebras}},
Math. Ann.
  \textbf{233} (1978), 145--153.

\bibitem{Dad:cont-fields}
M.~D\u{a}d\u{a}rlat, \emph{{Continuous fields of $C^{*}$-algebras
over
  finite-dimensional spaces}}, Preprint, 2006.

\bibitem{DadWinter:asu}
M.~D\u{a}d\u{a}rlat and W.~Winter, \emph{{On the $KK$-theory of
strongly
  self-absorbing $C^{*}$-algebras}}, in preparation, 2006.

\bibitem{GongJiangSu:Z}
G.~Gong, X.~Jiang, and H.~Su, \emph{{Obstructions to
$\mathcal{Z}$-stability
  for unital simple $C^*$-algebras}}, Canadian Math. Bull. \textbf{43} (2000),
  no.~4, 418--426.

\bibitem{GooHan:extending}
K.~R. Goodearl and D.~Handelman, \emph{{Rank functions and $K_0$
of regular
  rings}}, J. Pure Appl. Algebra \textbf{7} (1976), 195--216.

\bibitem{Haa:quasi}
U.~Haagerup, \emph{{Every quasi-trace on an exact $C^*$-algebra is
a trace}},
  preprint, 1991.

\bibitem{HirshbergWinter:Rokhlin-ssa}
I.~Hirshberg and W.~Winter, \emph{{Rokhlin actions and
self-absorbing
  $C^{*}$-algebras}}, Preprint, Math. Archive math.OA/0505356, 2005.

\bibitem{HjeRor:stable}
J.~Hjelmborg and M.~R{\o}rdam, \emph{{On stability of
$C^*$-algebras}}, J.
  Funct. Anal. \textbf{155} (1998), no.~1, 153--170.

\bibitem{HurWall:Dim}
W.~Hurewicz and H.~Wallmann, \emph{{Dimension Theory}}, Princeton
University
  Press, Princeton, 1941.

\bibitem{Hus:fibre}
D.~Husemoller, \emph{{Fibre Bundles}}, 3rd. ed., Graduate Texts in
Mathematics,
  no.~20, Springer Verlag, New York, 1966, 1994.

\bibitem{JiaSu:Z}
X.~Jiang and H.~Su, \emph{{On a simple unital projectionless
$C^*$-algebra}},
  American J. Math. \textbf{121} (1999), no.~2, 359--413.

\bibitem{Kas:Novikov}
G.~G. Kasparov, \emph{{Equivariant $KK$-theory and the Novikov
conjecture}},
  Invent. Math. \textbf{91} (1988), 147--201.

\bibitem{Kir:CentralSequences}
E.~Kirchberg, \emph{{Central sequences in $C^*$-algebras and
purely infinite
  $C^*$-algebras}}, Preprint.

\bibitem{Kir:spi}
\bysame, \emph{{Permanence properties of purely infinite
$C^*$-algebras}}, in
  preparation.

\bibitem{Kir:fields}
\bysame, \emph{{The classification of Purely Infinite
$C^*$-algebras using
  Kasparov's Theory}}, in preparation.

\bibitem{KirPhi:classI}
E.~Kirchberg and N.~C. Phillips, \emph{{Embedding of exact
$C^*$-algebras into
  ${\mathcal{O}}_2$}}, J. Reine Angew. Math. \textbf{525} (2000), 17--53.

\bibitem{KirRor:pi}
E.~Kirchberg and M.~R{\o}rdam, \emph{{Non-simple purely infinite
  $C^*$-algebras}}, American J. Math. \textbf{122} (2000), 637--666.

\bibitem{Phi:class}
N.~C. Phillips, \emph{{A Classification Theorem for Nuclear Purely
Infinite
  Simple $C^*$-Algebras}}, Documenta Math. \textbf{5} (2000), 49--114.

\bibitem{Ror:UHF}
M.~R{\o}rdam, \emph{{On the Structure of Simple $C^*$-algebras
Tensored with a
  UHF-Algebra}}, J. Funct. Anal. \textbf{100} (1991), 1--17.

\bibitem{Ror:UHFII}
\bysame, \emph{{On the Structure of Simple $C^*$-algebras Tensored
with a
  UHF-Algebra, II}}, J. Funct. Anal. \textbf{107} (1992), 255--269.

\bibitem{Ror:sns}
\bysame, \emph{{Stability of $C^*$-algebras is not a stable
property}},
  Documenta Math. \textbf{2} (1997), 375--386.

\bibitem{Ror:encyc}
\bysame, \emph{{Classification of Nuclear, Simple
$C^*$-algebras}},
  {Classification of Nuclear $C^*$-Algebras. Entropy in Operator Algebras}
  (J.~Cuntz and V.~Jones, eds.), vol. 126, Encyclopaedia of Mathematical
  Sciences. Subseries: Operator Algebras and Non-commutative Geometry, no. VII,
  Springer Verlag, Berlin, Heidelberg, 2001, pp.~1--145.

\bibitem{Ror:extension}
\bysame, \emph{{Extensions of stable $C^*$-algebras}}, Documenta
Math
  \textbf{6} (2001), 241--246.

\bibitem{Ror:simple}
\bysame, \emph{{A simple $C^*$-algebra with a finite and an
infinite
  projection}}, Acta Math. \textbf{191} (2003), 109--142.

\bibitem{Ror:stable}
\bysame, \emph{{Stable $C^*$-algebras}}, Advanced Studies in Pure
Mathematics
  38 ``Operator Algebras and Applications", Math. Soc. Japan, Tokyo, 2004,
  pp.~177--199.

\bibitem{Ror:Z-absorbing}
\bysame, \emph{{The stable and the real rank of
$\mathcal{Z}$-absorbing
  C*-algebras}}, International J.\ Math. \textbf{15} (2004), no.~10,
  1065--1084.

\bibitem{TomsWinter:ssa}
A.~Toms and W.~Winter, \emph{{Strongly self-absorbing
C$^*$-algebras}},
  Preprint, Math. Archive math.OA/0502211, to appear in Trans. Amer. Math.
  Soc., 2005.

\bibitem{TomsWinter:Zash}
\bysame, \emph{{$\Zh$-stable ASH $C^*$-algebras}}, Preprint, Math.
Archive
  math.OA/0508218, to appear in Can. J. Math., 2005.

\bibitem{TomsWinter:VI}
\bysame, \emph{{The Elliott conjecture for Villadsen algebras of
the first
  type}}, in preparation, 2006.

\bibitem{Winter:Z-class}
W.~Winter, \emph{{On the classification of simple
$\mathcal{Z}$-stable
  $C^*$-algebras with real rank zero and finite decomposition rank}}, J. London
  Math. Soc. \textbf{74} (2006), 167--183.

\bibitem{Winter:lfdr}
\bysame, \emph{{Simple $C^*$-algebras with locally finite
decomposition rank}},
  Preprint, Math. Archive math.OA/0602617, 2006.

\end{thebibliography}
\providecommand{\bysame}{\leavevmode\hbox
to3em{\hrulefill}\thinspace}
\providecommand{\MR}{\relax\ifhmode\unskip\space\fi MR }
\providecommand{\MRhref}[2]{%
  \href{http://www.ams.org/mathscinet-getitem?mr=#1}{#2}
} \providecommand{\href}[2]{#2}

\end{document}